\documentclass[11pt]{article}
\usepackage{amssymb}
\usepackage{amsmath}
\usepackage{pst-all}
\usepackage{graphicx}
\usepackage{pstricks}
\usepackage{makeidx}
\usepackage{enumerate}
\usepackage{mathtools}
\usepackage{comment}
\usepackage{hyperref}
\usepackage[a4paper,left=2cm,right=2cm,top=2.5cm,bottom=2.5cm]{geometry}

\newtheorem{lemma}{Lemma}[section]
\newtheorem{theorem}[lemma]{Theorem}
\newtheorem{corollary}[lemma]{Corollary}
\newtheorem{proposition}[lemma]{Proposition}
\newtheorem{remark}[lemma]{Remark}
\newtheorem{definition}[lemma]{Definition}
\newtheorem{example}[lemma]{Example}
\newtheorem{problem}[lemma]{Problem}

\def\bnum{\begin{enumerate} }
\def\enum{\end{enumerate}}
\def\bdf{\begin{definition}\rm }
\def\edf{\end{definition}}
\def\br{\begin{remark}\rm }
\def\er{\end{remark}}
\def\be{\begin{equation}}
\def\ee{\end{equation}}
\def\bt{\begin{theorem}}
\def\et{\end{theorem}}
\def\bl{\begin{lemma}}
\def\el{\end{lemma}}
\def\bc{\begin{corollary}}
\def\ec{\end{corollary}}
\def\bp{\begin{proposition}}
\def\ep{\end{proposition}}
\def\bxa{\begin{example}\rm }
\def\exa{\end{example}}
\def\ba{\begin{array}}
\def\ea{\end{array}}
\def\ben{\begin{eqnarray*}}
\def\een{\end{eqnarray*}}
\def\bdsc{\begin{description}}
\def\edsc{\end{description}}
\def\bpsp{\begin{pspicture}}
\def\epsp{\end{pspicture}}
\def\bea{\begin{eqnarray}}
\def\eea{\end{eqnarray}}
\def\btab{\begin{tabular}}
\def\etab{\end{tabular}}
\def\bpm{\begin{problem}}
\def\epm{\end{problem}}
\def\bfig{\begin{figure}}
\def\efig{\end{figure}}

\def\bnum{\begin{enumerate}\itemsep=0cm}
\def\enum{\end{enumerate}}

\def\pr{{\em Proof. }}

\def \qe{\hfill \vrule height4pt width 4pt depth 0pt}

\def\1{1\!\hspace{-.08cm}1}

\title{On the signless Laplacian spectrum of $k$-uniform hypergraphs}

\author{R.~B.~Bapat\thanks{Stat and Math Unit, ISI Delhi, New Delhi-110016, India}
\and
S.~S.~Saha\thanks{Department of Mathematics,
IIT Kharagpur, Kharagpur-721302, India.}
\and
S.~K.~Panda\thanks{Corresponding author: Department of Mathematics,
IIT Kharagpur, Kharagpur-721302, India. (spanda@maths.iitkgp.ac.in)}
}

\begin{document}
\pagestyle{myheadings}
\markboth{S.~K.~Panda}{On the signless Laplacian spectrum of $k$-uniform hypergraphs}
\maketitle
\begin{abstract}
Let $\mathcal{H}$ be a connected $k$-uniform hypergraph on $n$ vertices and $m$ hyperedges. In [A.~Banerjee, On the spectrum of hypergraph, Linear Algebra and its Application, 614(2021), 82--110], Anirban Banerjee introduced a new adjacency matrix for hypergraphs. In this article we consider the corresponding signless Laplacian matrix $Q(\mathcal{H})$. 
Let $q_{\max}=q_1\geq\cdots\geq q_n=q_{\min}$ be the eigenvalues of $Q(\mathcal{H})$ and let $d_{\max}=d_1\geq\cdots\geq d_n=d_{\min}$ be the degrees of $\mathcal{H}.$ We obtain some upper and lower bounds for $q_{\max}$ for connected $k$-uniform hypergraphs. The following are some of them. We prove that  $q_{\max}(\mathcal{H})\geq \frac{d_{\min}+\sqrt{d^{2}_{\min}+\frac{8}{k-1}T_{\min}}}{2}$
and we characterize the equality case. We prove that  $q_{\max}(\mathcal{H})\geq\frac{1}{k-1}\big((n-\tau)\binom{n-\tau}{k-1}+\tau\binom{n-\tau-1}{k-2}\big)$, where $\tau$ is weak independence number and we characterize the equality case. We prove that  $q_{\max}(\mathcal{H})\geq\frac{km}{n}\big(1+\frac{1}{\chi-1}\big)$ for a $k$-uniform hypergraph, where $\chi(\mathcal{H})$ is the strong chromatic number. We prove that $q_{\max}(\mathcal{H})\leq \max\limits_{i\in V}(d_i+\frac{s_i}{k-1}),$ where $s_i=\sum\limits_{j\sim i} \frac{d_jd_{ij}}{d_i}$ and $d_{ij}$ number of hyperedges containing $i$ and $j$ and we characterize the equality case.

In [D.~Cvetkovic, P.~Rowlinson and S.~K.~Simic, Signless Laplacians of finite graphs, Linear Algebra and its Applications 423 (2007) 155–171], Cvetkovic et al. proved that the least eigenvalue of the signless Laplacian of a connected graph is equal to $0$ if and only if the graph is bipartite. We prove that connected bipartite graphs are the only $k$-uniform connected hypergraphs for which the least eigenvalue of the signless Laplacian is $0.$ Next we prove that $q_{\min}(\mathcal{H})\leq d_{\max}-\frac{1}{k-1}.$ We prove that  $q_{\min}(\mathcal{H})\leq2\sqrt{\frac{Z_{1}(\mathcal{H})}{n}},$ where $Z_{1}(\mathcal{H})$ is first Zagreb index.
We also obtain some bounds for the least eigenvalue of $k$-uniform hypergraph. 

Carla Silva Oliveira et al.[Carla Silva Oliveira,  Leonardo Silva de Lima, Nair Maria Maia de Abreu and Steve Kirkland,Linear Algebra and its Applications 432(2010), 2342--2351] proved that, for any graph on $n\geq 2$ vertices with at least one edge, $s_Q(G)\geq 2$ and characterized the equality case. In this article we prove that, for any connected $k$-uniform hypergraph $\mathcal{H}$ on $n\geq 2$ vertices with at least one hyperedge, $s_Q(\mathcal{H})>1$. We prove that $\inf\big\{s_Q(\mathcal{H}): \mathcal{H}\in \mathcal{F}\big\}=1,$ where $\mathcal{F}$ is the set of all hypergraphs with finite vertices. We prove that $s_Q(\mathcal{H})\geq \frac{2nd_{\min}-km}{n-1}$ and we characterize the equality case.
We obtain a lower bound for $s_Q(\mathcal{H})$ in terms of frist Zagreb index and characterize the equality case. We obtain an upper bound and a lower bound for the $s_Q(\mathcal{H})$  for a $k$-uniform hypergraph in terms of strong chromatic number. 
\end{abstract}

\sloppy

\section{Introduction}
In 1973, Berge Claude introduced the concept of Hypergraph \cite{claude1973}. An {\em undirected hypergraph} $\mathcal{H}$ is a pair $\mathcal{H}= (V,E)$, where $V$ is the set of $n$ elements called the vertex set and $E$ is a set of non-empty subsets of $V$ is called the hyperedge set. The {\em sub-hypergraph} induced by the set of vertices $V'$ is defined as the pair $(V', E')$, where $E'$ is a subset of E and $V'$ is a nonempty set of vertices of the hypergraph. A hypergraph is {\em finite} if $V$ is finite set and with no vertices is called a {\em null hypergraph}, with only one vertex is called a {\em trivial hypergraph}, and all other hypergraphs are nontrivial. A hypergraph with no hyperedges is called {\em empty}. A hypergraph is said to be {\em simple} if for any $e_i,e_j\in E$ such that $e_i\neq e_j$ and nether $e_i\subset e_j$ nor $e_j\subset e_i$.

For $i,j\in V$ , $i$ and $j$ are said to be {\em adjacent} if there exists an hyperedge $e\in E$ such that $i,j\in e$. We use $i\sim j$ to denote the vertex $i$ is adjacent to the vertex $j$. 
Two hyperedges are said to be {\em adjacent} if their intersection is not empty. A vertex $i$ and an hyperedge $e_i$ are said to be {\em incident} if $i\in e_i$. Otherwise, the vertex is non-incident to the hyperedge. The size of an hyperedge $e\in E$ is refer as the number of vertices in the hyperedge $e$. A hypergraph $\mathcal{H}$ is called a {\em $k$-uniform} hypergraph for any integer $k\geq 2$, if for all $e\in E$,  $|e|=k$. A $k$-uniform hypergraph $\mathcal{H}$ with n vertices is said to be {\em complete} if $\mathcal{H}$ contains all the hyperedges of size $k$ and it is denoted by $K^k_n.$

For vertices $i$ and $j$ (not necessarily distinct) in hypergraph $\mathcal{H}$, $(i,j)$-{\em walk} of length $k$ in $\mathcal{H}$ is an alternating sequence $i= i_0~e_1~i_1~e_2~i_2\ldots e_k~i_k= j$ of vertices and hyperedges, starting at $i$ and ending at $j$ such that every hyperedge is incident on vertices preceding and following it, that is, $v_i,v_{i-1}\in e_i$, $1\leq i\leq k$. If $i=j$ then the walk is called a {\em closed walk}, otherwise an open walk. A path is an {\em open walk} where all the vertices and hyperedges are distinct. A {\em cycle} is a closed walk where all the vertices are distinct except the end vertices. A hypergraph $\mathcal{H}$ is said to be {\em connected} if for any two distinct vertices $i,j$ there exists $(i,j)$-path.  Throughout this article we only consider undirected, finite, simple, non-empty, nontrivial and connected hypergraphs.

Let $E_{ij}=\{e\in E: i,j\in e\}$. The {\em co-degree} $d_{ij}$ of vertices $i$ and $j$ is the cardinality of $E_{ij}$. The degree of a vertex is the number of hyperedges that contain the vertex.  Let $d_{\max}=d_1\geq d_2\geq\cdots\geq d_n=d_{\min}$ be the dgeree sequence of $\mathcal{H}.$ A hypergraph $\mathcal{H}(V,E)$ is called {\em bipartite} if $V$ can be partitioned into two disjoint subsets $V_1$ and $V_2$ such that for each hyperedge $e\in E$, $e\cap V_1\neq \emptyset$ and $e\cap V_2\neq \emptyset$. An {\em $m$-uniform complete bipartite} hypergraph is denoted by $C^m_{m_1,m_2}$ , where $|V_1|=m_1$ and $|V_2|=m_2$.

\begin{definition}
The $2$-degree of a vertex $i$ in a hypergraph $\mathcal{H}$ is the sum of the degrees of the vertices adjacent to $i$ and denoted by $T_{i}=\sum\limits_{j\sim i}d_{j}$. We use $T_{\min}$ to denote $\min\limits_{i\in V}T_{i}.$
\end{definition}

\begin{definition}
Let $\mathcal{H}=(V,E)$ be a hypergraph. A subset $S$ of the vertex set $V$ is called an independet set/coclique if for every two vertices in $S$, there is no hyperedge in $\mathcal{H}$ containing them. The {\em independence number} of $\mathcal{H}$ is the maximum cardinality of an independent set of vertices in $\mathcal{H}.$
An independent set $S$ is called {\em weak independent set} if each vertex of $S$ is adjacent to each vertex of $V\setminus S.$  The {\em weak independence number} of a hypergraph $\mathcal{H}$ is the maximum cardinality of a weak independent set of vertices in $\mathcal{H}.$
and we denote it by $\tau(\mathcal{H})$.
\end{definition}

\begin{definition}
The strong chromatic number $\chi(\mathcal{H})$ of a hypergraph $\mathcal{H}$ is the minimum number of colors required to color the vertices so that adjacent vertices get distinct color.
\end{definition}

\begin{definition}
Let $i$ be a vertex of a hypergraph $\mathcal{H}$. Then the {average degree of the vertex} $i$ is $\sum\limits_{j\sim i}\frac{d_jd_{ij}}{d_i}$ and it is denoted by $s_i$.
\end{definition}

In 2021, Anirban Banerjee \cite{banerjee2021} was first introduced a new type of adjacency matrix $\mathcal{A}(\mathcal{H})$ and corresponding Laplacian matrix $\mathcal{L}(\mathcal{H})$ for hypergraphs. 
The \textbf {adjacency matrix} $\mathcal{A}({\mathcal{H}})= [(\mathcal{A}({\mathcal{H}}))_{ij}]$ of a hypergraph $\mathcal{H}= (V,E)$ is defined as 

$$
(\mathcal{A}({\mathcal{H}}))_{ij}=\begin{cases}
    \sum\limits_{e\in E_{ij}}\frac{1}{|e|-1}~~~,\mbox{ if }~~i\sim j,\\
     0~~~~~~~~~~~~~~~,\mbox{ otherwise.}
 \end{cases}$$

Then $d_i=\sum\limits_{i\sim j}(\mathcal{A}({\mathcal{H}}))_{ij}$ is the degree of the vertex $i$.
The Laplacian matrix of $\mathcal{H}$ is defined as  $\mathcal{L}(\mathcal{H})=\mathcal{D}(\mathcal{H})-\mathcal{A}(\mathcal{H})$, where $\mathcal{D}(\mathcal{H})=diag\big( d_1(\mathcal{H}),\ldots, d_n(\mathcal{H})\big)$. The matrix $\mathcal{L}$ is symmetric and positive semi-definite. For more details, see \cite{banerjee2021,saha2022}. 

In this article we consider the corresponding {\em signless Laplacian matrix}. 
The \textbf {signless Laplacian matrix} of a hypergraph $\mathcal{H}$ is denoted by $Q(\mathcal{H})=[q({\mathcal{H}})_{ij}]_{n\times n}$ and is define by $Q(\mathcal{H})=\mathcal{D}(\mathcal{H})+\mathcal{A}(\mathcal{H})$, 
The signless Laplacian matrix is a symmetric and positive semi-definite matrix. The spectrum of the signless Laplacian matrix of graphs is well studied by several authors \cite{Andrade2019,das2004,das2010,liu2010,lima2011,Liu2010,MadenKinkar,Oliveira2010,Zhang2009}. Let $q_{\max}=q_1\geq q_2\geq \cdots\geq q_n=q_{\min}$ be the eigenvalues of $Q(\mathcal{H}).$

We arrange the article as follows. 
In Section 2, we prove that  $q_{\max}(\mathcal{H})\geq \frac{d_{\min}+\sqrt{d^{2}_{\min}+\frac{8}{k-1}T_{\min}}}{2}$
and we characterize the equality case. We prove that  $q_{\max}(\mathcal{H})\geq\frac{1}{k-1}\big((n-\tau)\binom{n-\tau}{k-1}+\tau\binom{n-\tau-1}{k-2}\big)$, where $\tau$ is weak independence number and we characterize the equality case. For a connected $k$-uniform hypergraphs prove that  $q_{\max}(\mathcal{H})\geq\frac{km}{n}\big(1+\frac{1}{\chi-1}\big)$, where $\chi(\mathcal{H})$ is the strong chromatic number. We prove that $q_{\max}(\mathcal{H})\leq \max\limits_{i\in V}(d_i+\frac{s_i}{k-1}),$ where $s_i$ is the average degree of the vertex $i$ and we characterize the equality case. We also obtain some other bounds for the sigless Laplacian spectral radius. 
In Section 3, we prove that connected bipartite graphs are the only $k$-uniform connected hypergraphs for which the least eigenvalue of the signless Laplacian is $0.$ Next we prove that $q_{\min}(\mathcal{H})\leq d_{\max}-\frac{1}{k-1}.$ We prove that  $q_{\min}(\mathcal{H})\leq2\sqrt{\frac{Z_{1}(\mathcal{H})}{n}},$ where $Z_{1}(\mathcal{H})$ is first Zagreb index.
We also obtain some other bounds for the least eigenvalue of $k$-uniform hypergraphs. 
In Section 4, we prove that, for any $k$-uniform hypergraph $\mathcal{H}$ on $n\geq 2$ vertices with at least one hyperedge, $s_Q(\mathcal{H})>1$.  We prove that $\inf\big\{s_Q(\mathcal{H}): \mathcal{H}\in \mathcal{F}\big\}=1,$ where $\mathcal{F}$ is the set of all hypergraphs with finite vertices.
We prove that $s_Q(\mathcal{H})\geq \frac{2nd_{\min}-km}{n-1}$ and we characterize the equality case.
We obtain a lower bound for $s_Q(\mathcal{H})$ in terms of frist Zagreb index and all the $k$-uniform hypergraphs of order $n$ which attain this lower bound are characterized. We obtain an upper bound and a lower bound for the $s_Q(\mathcal{H})$  for a $k$-uniform hypergraph in terms of strong chromatic number.

\section{Largest signless Laplacian eigenvalue of a $k$-uniform hypergraph}
In this section, we shall try to obtain bounds for the signless Laplacian spectral radius of $k$-uniform hypergraph.
\begin{definition}
For $1\leq k\leq n$ denote by $E(k)$ the set of
all $k$-elementary subsets of $V$. Let $\mathcal{H}=(V(\mathcal{H}),E(\mathcal{H}))$. The complement $\bar{\mathcal{H}}$ of the hypergraph $\mathcal{H}$ is defined as the pair $\bigg(V(\mathcal{H}),E(k)\setminus E(\mathcal{H})\bigg)$.
\end{definition}

\begin{remark} \label{rem1}
Let $\cal{H}$ be a $k$- uniform hypergraph. Then the following are true.

\begin{enumerate}
\item The complement hypergraph $\mathcal{\overline{H}}$ is also $k$-uniform hypergraph.

\item Let $d_{i}$ be the degree of $i$ in $\mathcal{H}$ and $\overline{d_{i}}$ be the degree of $i$ in $\mathcal{\overline{H}}$. Then $d_{i}+\overline{d_{i}}=\theta_{1} =\left( \begin{array}{c} n-1 \\ k-1 \end{array} \right)$. Also $\mathcal{A}({\mathcal{H}})+\mathcal{A}({\mathcal{\overline{H}}})=\theta (J_{n}-I_{n})$, where $\theta=\frac{\left( \begin{array}{c} n-2 \\ k-2 \end{array} \right)}{k-1}$ and $J_n$ is a matrix of size $n$ whose all entries are $1$.
\end{enumerate}
\end{remark}

\begin{lemma}{\rm\cite{horn2013}}\label{thm2.3}
Let $A,B$ be two Hermitian matrices. Let $\lambda_{1}(A)\geq\lambda_{2}(A)\geq\cdots\geq\lambda_{n}(A)$, $\lambda_{1}(B)\geq\lambda_{2}(B)\geq\cdots\geq\lambda_{n}(B)$ and  $\lambda_{1}(A+B)\geq\lambda_{2}(A+B)\geq\cdots\geq\lambda_{n}(A+B)$ be the eigenvalues of $A$, $B$ and $A+B$, respectively. Then
\begin{center}
    $\lambda_{i}(A)+\lambda_{n}(B)\leq\lambda_{i}(A+B)\leq\lambda_{i}(A)+\lambda_{1}(B)$, $i=1,2,\ldots,n.$
\end{center}

\end{lemma}

The following theorem gives a relation between the spectral radius of signless Laplacian and smallest signless Laplacian eigenvalues of a $k$-uniform hypergraphs $\mathcal{H}$ and $\overline{\mathcal{H}}$ respectively.

\begin{theorem}\label{thm:2.4}
Let $\mathcal{H}$ be a $k$-uniform hypergraph with n vertices and $\overline{\mathcal{H}}$ be the complement of $\mathcal{H}$. Let $q_{\max}(\overline{\mathcal{H}})=q_1(\overline{\mathcal{H}})\geq q_2(\overline{\mathcal{H}})\geq\cdots\geq q_n(\overline{\mathcal{H}}) =q_{\min}(\overline{\mathcal{H}})$ be the signless Laplacian eigenvalues of $\overline{\mathcal{H}}.$ Then $q_{\max}(\overline{\mathcal{H}})\geq(n-2)\theta-q_{\min}(\mathcal{H})$.
\end{theorem}

\pr Since we have, 
\begin{align*}
 Q(\mathcal{H})+Q({\mathcal{\overline{H}}})&=\mathcal{D}(\mathcal{H})+\mathcal{D}(\overline{\mathcal{H}})+\mathcal{A}(\mathcal{H})+\mathcal{A}(\overline{\mathcal{H}})\\
&=\theta_{1}I_{n}+\theta (J_{n}-I_{n}), \mbox{by Remark (\ref{rem1})}\\
&=(\theta_{1}-\theta)I_{n}+\theta J_{n}
=(n-2)\theta I_{n}+\theta J_{n}
\end{align*} 

Therefore $Q(\overline{\mathcal{H}})=(n-2)\theta I_{n}-Q(\mathcal{H})+\theta J_{n}$.
Let $A=(n-2)\theta I_{n}-Q(\mathcal{H})$ and $B=\theta J_{n}$.
 
 Using Lemma \ref{thm2.3}, $q_{\max}(\overline{\mathcal{H}})\geq(n-2)\theta-q_{\min}(\mathcal{H})+0.$ Hence $q_{\max}(\overline{\mathcal{H}})\geq(n-2)\theta-q_{\min}(\mathcal{H})$.\qe

The following three results will be used in the sequel.

\begin{lemma}{\rm\cite{vol2009}}\label{thm2.5}
 Let $\mathcal{H}$ be a $k$-uniform hypergraph with $m$ hyperedges. Then

\begin{center}
    $\sum\limits_{i \in V} d_i = k m.$ 
\end{center}
\end{lemma} 

\begin{theorem}{\rm\cite{saha2022}}\label{thm:2.6}
Let $\mathcal{H}=(V,E)$ be a $k$-uniform hypergraph.
Then 
\begin{center}
    $\mu_{\max}(\mathcal{H})\geq d_{\max}+\frac{1}{(k-1)},$
    
\end{center}
where $\mu_{\max}(\mathcal{H})$ and $d_{\max}$ are the Laplacian spectral radius and maximum degree of $\mathcal{H}$, respectively.
\end{theorem}

\begin{lemma}{\rm\cite{saha2022}}\label{lem2.7}
Let $\mathcal{H}$ be a hypergraph. Then $\mu_{\max}({\mathcal{H}})\leq q_{\max}(\mathcal{H})$,
where $\mu_{\max}(\mathcal{H})$ and $q_{\max}(\mathcal{H})$ are the spectral radius of the Laplacian and signless Laplacian matrix of $\mathcal{H},$ respectively.

\end{lemma}

To proceed further we need the following lemma.

\begin{lemma}\label{newlem2} {\rm\cite{ellingham2000}}\label{lem3.7}
Let $\cal{H}$ be a hypergraph with a signless Laplacian matrix $Q(\mathcal{H})$. Let $P$ be any polynomial and $R_{i}(P(Q(\mathcal{H}))$ be the $ith$-row sum of $P(Q(\mathcal{H}))$. Then
\begin{center}
    $\min\limits_{i\in V}R_{i}(P(Q(\mathcal{H}))\leq P(q_{\max}(\mathcal{H}))\leq\max\limits_{i\in V}R_{i}(P(Q(\mathcal{H})).$
\end{center}
\end{lemma}

The following theorem gives a lower bound for the signless Laplacian spectral radius of a $k$-uniform hypergraph in terms of minimum degree and 2-degree. We use $R_{i}(M)$ to denote the $i$th row sum of the matrix $M$.

\begin{theorem}\label{thm 2.15}
Let $\cal{H}$ be a connected $k$-uniform hypergraph. Let $d_{\min}$ be a minimun degree of $\cal{H}$. Then
\begin{center}
    $q_{\max}(\mathcal{H})\geq \frac{d_{\min}+\sqrt{d^{2}_{\min}+\frac{8}{k-1}T_{\min}}}{2},$
\end{center}
where $T_{\min}=\min\limits_{i\in V}T_{i}.$ Equality holds if and only if $\mathcal{H}$ is regular and for any two adjacent vertices there is exactly one hyperedge containing those two vertices.
\end{theorem}

\pr Let $i$ be any vertex of $\mathcal{H}$. Then $R_{i}(Q(\mathcal{H}))=2d_{i}$ and  $R_{i}(\mathcal{A}^{2}(\mathcal{H}))=R_{i}(\mathcal{A}(\mathcal{H})\mathcal{D}(\mathcal{H}))=\frac{1}{(k-1)}\sum\limits_{i\sim j}d_{ij}d_{j}.$ 

Now $Q^{2}(\mathcal{H})= \big(\mathcal{D}(\mathcal{H})+\mathcal{A}(\mathcal{H})\big)\big(\mathcal{D}(\mathcal{H})+\mathcal{A}(\mathcal{H})\big)$
$=\mathcal{D}(\mathcal{H})Q(\mathcal{H})+\mathcal{A}(\mathcal{H})\mathcal{D}(\mathcal{H})+\mathcal{A}^2(\mathcal{H})$ 

Taking $ith$-row sum both sides we get,
\begin{align*} 
R_i(Q^2(\mathcal{H}))&=R_i(\mathcal{D}(\mathcal{H})Q(\mathcal{H}))+2R_i(\mathcal{A}^2(\mathcal{H})).\\
&=d_{i}R_{i}\big(Q(\mathcal{H})\big)+\frac{2}{(k-1)}\sum\limits_{i\sim j}d_{ij}d_{j}\\
&\geq d_{\min} R_{i}(Q(\mathcal{H}))+\frac{2}{k-1}\sum\limits_{i\sim j}d_{j} \hspace{.5cm}(\mbox{as}\hspace{.2cm} d_{ij}\geq 1)\\
&\geq d_{\min}R_{i}(Q(\mathcal{H})+\frac{2}{k-1}T_{i}.
\end{align*}

So $R_i(Q^2(\mathcal{H}))\geq d_{\min}R_{i}(Q(\mathcal{H})+\frac{2}{k-1}T_{i}$ for $i=1,\ldots,n.$

This implies $R_{i}\Big(Q^2(\mathcal{H})-d_{\min}Q(\mathcal{H})\Big)\geq\frac{2}{k-1}T_{i}$ for $i=1,\ldots,n.$ 
Taking minimum over all the vertices we get, 
\begin{center}
$\min\limits_{i\in V}R_{i}\big(Q^2(\mathcal{H})-d_{\min}Q(\mathcal{H})\big)\geq\frac{2}{k-1}\min\limits_{i\in V}T_{i}.$
\end{center}

Using Lemma \ref{lem3.7} we have, $q_{\max}(\mathcal{H})^{2}-d_{\min} q_{\max}(\mathcal{H})-\frac{2}{k-1}T_{\min}\geq0.$ 
Hence it follows that, 
\begin{center}
$q_{\max}(\mathcal{H})\geq\frac{d_{\min}+\sqrt{d^{2}_{\min}+\frac{8}{k-1}T_{\min}}}{2}$.
\end{center}

Suppose equality holds. Then form the above proof we have $R_i(Q^2(\mathcal{H}))= d_{\min}R_{i}(Q(\mathcal{H})+\frac{2}{k-1}T_{i}$, $d_i=d_{\min}$ for $i=1,\ldots,n$ and $d_{ij}=1$ for $i,j=1,\ldots n$ and $i\neq j.$ This implies $\mathcal{H}$ is regular and for any two adjacent vertices there is exactly one hyperedge containing those two vertices. 

Conversely, if $\mathcal{H}$ is regular and for any two adjacent vertices there is exactly one hyperedge containing those two vertices. Then $d_1=d_2=\cdots=d_n=d$ and $T_{i}=\sum\limits_{i\sim j}d_{j}=(k-1)d^2$. Therefore $T_{\min}=(k-1)d^2$. The right hand side of the inequality is $2d$ and $q_{\max}=2d$. Hence the equality holds. 
\qe
 \begin{lemma}\rm\cite{saha2022} \label{l3.21}
Let $\mathcal{H}=(V,E)$ be a connected hypergraph on $n$ vertices and $E'$ be a subset of $E$. Then sub-hypergraph of $\mathcal{H}$ obtained by deleting all the hyperedges in $E'$ is denoted $\mathcal{H}-E'$. Then
\begin{center}
    $\mu_{i}(\mathcal{H})\geq\mu_{i}(\mathcal{H}-E')$, for $i=1,2,\ldots n$.
\end{center}
  \end{lemma}

  \begin{definition}\rm\cite{banerjee2021}
Let $\mathcal{H}=(V,E)$ be a $k$-uniform hypergraph and $V_{1},V_{2},\ldots,V_{t}$ be a partition of the vertex set $V$. Then $V_{1},V_{2},\ldots,V_{t}$ is said to form an {\em equitable partition} if for each $p,q\in \{1,2,\ldots, t\}$ and for any $i\in V_{p}$
\begin{center}
    $\sum\limits_{j;j\in V{q}}(\mathcal{A}({\mathcal{H}}))_{ij}=m_{pq},$
\end{center}
where $m_{pq}$ are constant depending on $p$ and $q$.
\end{definition}
  
  \begin{theorem}\rm\cite{saha2022}\label{t3.23}
 Let $V_{1},V_{2},\ldots,V_{t}$ be an equitable partition of the vertex set of a uniform hypergraph $\mathcal{H}$ with a constants $m_{pq}$ ($p,q=1,2,\ldots,t$) and let $M=[(m_{pq})]$ be the $t\times t$ matrix defined by

$$
m_{pq}=\begin{cases}
     -m_{pq}~~~~~~~~~~~~~~~~,\mbox{ if }~~~p\neq q,\\
     \big(\sum\limits_{s=1}^{t}m_{ps}\big)-m_{pp}~~,\mbox{ if}~~~~ p=q.
 \end{cases}$$

Then each eigenvalue of $M$ is also an eigenvalue of $\mathcal{L}({\mathcal{H}})$.
\end{theorem}
  
\begin{lemma}\rm\cite{horn2013}\label{l3.24}
Let $x,y\in\mathbb{C}^n$ and $a\in\mathbb{C}$. Consider the matrix $A = \begin{bmatrix} 
    \lambda I_n & x \\
     y^{*} & a \\
\end{bmatrix}$. Then the characteristic polynomial of A is $Ch_{A}(t)=(t-\lambda)^{n-1}(t^2-(a+\lambda)t+a\lambda-y^{*}x)$ .
\end{lemma}

Next we supply a lower bound of $Q$-spectral radius in terms of number of vertices and weak independence number.

\begin{theorem}
Let $\mathcal{H}=(V,E)$ be a connected $k$-uniform hypergraph with $n$ vertex. Let $\tau=\tau(\mathcal{H})$ be the weak independence number of $\mathcal{H}$. Then \begin{center}
    $q_{\max}(\mathcal{H})\geq\frac{1}{k-1}\bigg((n-\tau)\binom{n-\tau}{k-1}+\tau\binom{n-\tau-1}{k-2}\bigg),$
\end{center}
Equality holds if and only if $\mathcal{H}$ is a bipartite complete graph $K_{\tau, n-\tau}.$
\end{theorem}

\pr Let us take the vertex set $V=\{v_1,v_2,\ldots,v_n\}$ such that $v_1$,$v_2$,$\ldots$,$v_{\tau}$ are the vertices corresponding to  the maximum weak independent set of $\mathcal{H}$. 

Let $\mathcal{H}'$ be the sub-hypergraph of $\mathcal{H}$ such that remove all the hyperedges between the vertices $v_{\tau+1},v_{\tau+2},\ldots,v_{n}$.

Now taking $V_1=\{v_1\}$, $V_2=\{v_2\}$,$\ldots$,$V_{\tau}=\{v_{\tau}\}$ and $V_{\tau+1}=\{v_{\tau+1},v_{\tau+2},\ldots,v_{n}\}$. Then $V_{1},V_{2},\ldots,V_{\tau}$ and $V_{\tau+1}$ form an equitable partition of $\mathcal{H}'$. Then the matrix 

     $$M(\mathcal{H}') = \begin{bmatrix} 
    \frac{(n-\tau)\binom{n-\tau}{k-1}}{k-1} & 0 & 0 & \cdots & 0 & -\frac{(n-\tau)\binom{n-\tau}{k-1}}{k-1}\\
0  & \frac{(n-\tau)\binom{n-\tau}{k-1}}{k-1} & 0 & \cdots & 0 & -\frac{(n-\tau)\binom{n-\tau}{k-1}}{k-1}\\
\vdots & \vdots & \vdots & \cdots & \vdots & \vdots\\
0  & 0 & 0 & \cdots & \frac{(n-\tau)\binom{n-\tau}{k-1}}{k-1} & -\frac{(n-\tau)\binom{n-\tau}{k-1}}{k-1}\\
-\frac{\binom{n-\tau-1}{k-2}}{k-1} & -\frac{\binom{n-\tau-1}{k-2}}{k-1} & -\frac{\binom{n-\tau-1}{k-2}}{k-1} & \cdots & -\frac{\binom{n-\tau-1}{k-2}}{k-1} & \frac{\tau\binom{n-\tau-1}{k-2}}{k-1}\\
    \end{bmatrix}_{(\tau+1)\times(\tau+1)}$$
    
    $$\hspace{-4.7cm}= \begin{bmatrix} 
    \frac{(n-\tau)\binom{n-\tau}{k-1}}{k-1}I_{\tau} & -\frac{(n-\tau)\binom{n-\tau}{k-1}}{k-1}e_{\tau}  \\
    -\frac{\binom{n-\tau-1}{k-2}}{k-1}e^t_{\tau} & \frac{\tau\binom{n-\tau-1}{k-2}}{k-1} \\
\end{bmatrix}_{(\tau+1)\times(\tau+1)} $$
\hspace{1cm} where $e_{\tau}=(1,1,\ldots,1)^t$ and $I_{\tau}$ is the $\tau\times\tau$ identity matrix.  

Therefore by Lemma $\ref{l3.24}$, the characteristic polynomial is $Ch_{M}(t)=\bigg(t -\frac{(n-\tau)\binom{n-\tau}{k-1}}{k-1}\bigg)^{\tau-1}\bigg(t^2-\big(\frac{(n-\tau)\binom{n-\tau}{k-1}}{k-1}+\frac{\tau\binom{n-\tau-1}{k-2}}{k-1}\big)t\bigg)$. The eigenvalues of $M(\mathcal{H}')$ are 0,  $\bigg(\frac{(n-\tau)\binom{n-\tau}{k-1}}{k-1}+\frac{\tau\binom{n-\tau-1}{k-2}}{k-1}\bigg)$ with multiplicity 1 and  $\frac{(n-\tau)\binom{n-\tau}{k-1}}{k-1}$ with  multiplicity $\tau-1$, respectively. 

By Lemma $\ref{t3.23}$, if $\mu\in spec(M(\mathcal{H}'))$ implies $\mu\in spec(\mathcal{L}(\mathcal{H}'))$ and using Lemma  $\ref{l3.21}$, $\mu_{\max}(\mathcal{H})\geq\mu_{\max}(\mathcal{H}')=\bigg(\frac{(n-\tau)\binom{n-\tau}{k-1}}{k-1}+\frac{\tau\binom{n-\tau-1}{k-2}}{k-1}\bigg)$.

Hence by Lemma $\ref{lem2.7}$ we get, $q_{\max}(\mathcal{H})\geq\frac{1}{k-1}\bigg((n-\tau)\binom{n-\tau}{k-1}+\tau\binom{n-\tau-1}{k-2}\bigg)$.

We now prove the equality part. Suppose equality holds. Then from the above proof we have $q_{\max}(\mathcal{H})=\mu_{\max}(\mathcal{H})$ and $\mu_{\max}(\mathcal{H})=\mu_{\max}(\mathcal{H'}).$  Let $G$ be the weighted graph such that $Q(\mathcal{H})$ is signless Laplacian matrix. For this weighted graph we have $q_{\max}(G)=\mu_{\max}(G).$ Using the proof of \cite[Lemma 2.1]{zhang2002}, we have $G$ is bipartite. That is, $Q(\mathcal{H})=\begin{bmatrix} O & B^t\\ B & O\end{bmatrix}.$ This is not possible if $k>3.$ Hence $\mathcal{H}$ is bipartite graph.
The second equality says that there is no hyperedge between the vertices $V\backslash S.$ This implies that $\mathcal{H}$ is a complete bipartite graph $K_{\tau,n-\tau}.$

Converse part. If $\mathcal{H}$ is a complete bipartite graph $K_{\tau,n-\tau}$, then we can easily prove that $q_{\max}(\mathcal{H})=\frac{1}{k-1}\bigg((n-\tau)\binom{n-\tau}{k-1}+\tau\binom{n-\tau-1}{k-2}\bigg).$ \qe

\begin{theorem}
Let $\mathcal{H}$ be a connected $k$-uniform hypergraph with $n(\geq 2)$ vertices and $m$ hyperedges. Let $\chi=\chi(\mathcal{H})$ be the strong chromatic number of $\mathcal{H}.$ Then 
\begin{center}
    
    $q_{\max}(\mathcal{H})\geq\frac{km}{n} \big(1+\frac{1}{\chi-1}\big)$.
\end{center}
\end{theorem}

\pr We need $\chi=\chi(\mathcal{H})$ number of colors to color the vertices. Let $V_1,\ldots,V_{\chi}$ be the color partition. 
Since $Q(\mathcal{H})$ is symmetric, we have
$q_{\max}(\mathcal{H})=\max\limits_{x\in \mathbb{R}^n} \frac{x^tQ(\mathcal{H})x}{||x||^2}.$ Therefore 
$q_{\max}(\mathcal{H})||x||^2\geq x^tQ(\mathcal{H})x.$ Consider the vector $x=(x_1,\ldots,x_n)^t$ such that

$$
x(i)=\begin{cases}
     \chi -1~~~,\mbox{ if }~~~i\in V_k \\
     -1~~~~~~,\mbox{ elsewhere.}
 \end{cases}$$

Then 
\begin{align*}
||x||^2 &=(\chi-1)^2|V_k|+|V\setminus V_k|\\
        &=\chi(\chi-2)|V_k|+n.
\end{align*}

We now calculate $x^tQ(\mathcal{H})x=\sum\limits_{i\sim j}q_{ij} (x_i+x_j)^2.$

\begin{align*}
\sum\limits_{i\sim j}q_{ij} (x_i+x_j)^2&=\sum\limits_{\underset{i\in V_k, j\notin V_k}{i\sim j}}q_{ij} (x_i+x_j)^2+\sum\limits_{\underset{i\notin V_k, j\notin V_k}{i\sim j }}q_{ij} (x_i+x_j)^2\\
&=\sum\limits_{\underset{i\in V_k, j\notin V_k}{i\sim j}}q_{ij} (\chi^2-4\chi+4)+\sum\limits_{\underset{i\notin V_k, j\notin V_k}{i\sim j}}4q_{ij}\\
&=(\chi^2-4\chi)\sum\limits_{\underset{i\in V_k, j\notin V_k}{i\sim j}}q_{ij}+4\sum\limits_{\substack{i\sim j}}q_{ij}\\
&=(\chi^2-4\chi)e(V_k,V\setminus V_k)+4km, \big( \mbox{ using Lemma \ref{thm2.5} }\big)\\
\end{align*} 

where $e(V_k,V\backslash V_{k})=\sum\limits_{\underset{i\in V_k, j\notin V_k}{i\sim j}}q_{ij}$.
 
Therefore $q_{\max}(\mathcal{H})(\chi(\chi-2)|V_k|+n)\geq (\chi^2-4\chi)e(V_k,V\setminus V_k)+4km$, for all $k=1,\ldots,\chi.$
After summing them, we have
\begin{align*}
&q_{\max}(\mathcal{H})(\chi(\chi-2)n+\chi n)\geq \sum\limits_{k=1}^{\chi}(\chi^2-4\chi)e(V_k,V\setminus V_k)+4\chi km\\
&\Rightarrow q_{\max}(\mathcal{H})\chi(\chi-1)n\geq \chi^2km\\
&\therefore q_{\max}(\mathcal{H})\geq\frac{km}{n}\big(1+\frac{1}{\chi-1}\big).
\end{align*}
Hence the proof.
\qe

\begin{theorem}
Let $\mathcal{H}$ be a connected $k$-uniform hypergraph. Then 
\begin{center}
$q_{\max}(\mathcal{H})\leq \max\limits_{i\in V}(d_i+\frac{s_i}{k-1}),$ where $s_i=\sum\limits_{j\sim i} \frac{d_jd_{ij}}{d_i}.$
\end{center}
Equality holds if and only if either $\mathcal{H}$ is a regular hypergraph or a bipartite semiregular graph.
\end{theorem}

\pr Consider the matrix $\mathcal{D}(\mathcal{H})^{-1}Q(\mathcal{H})\mathcal{D}(\mathcal{H}).$ The $ij^{th}$ entry of $\mathcal{D}(\mathcal{H})^{-1}Q(\mathcal{H})\mathcal{D}(\mathcal{H})$ is

$$
\begin{cases}
     d_{i}~~~~~~~~,\mbox{ if }~~~i= j ,\\
     \frac{d_j}{d_i}\frac{d_{ij}}{k-1}~~~~,\mbox{if}~~~i\sim j\\
     0~~~~~~~~~~,\mbox{ elsewhere.}
 \end{cases}$$

Let $x=(x_1,\ldots,x_n)^t$ be an eigenvector corresponding to the eigenvalue $q_{\max}$ of $\mathcal{D}(\mathcal{H})^{-1}Q(\mathcal{H})\mathcal{D}(\mathcal{H}).$ Let $x_i=\max\limits_{j}x_j$. Set the vector $y=(y_1,\ldots,y_n)^t$ such that $y_k=\frac{x_k}{x_i}$ for $k=1,\ldots,n.$ Then $y$ is an eigenvector corresponding to the eigenvalue $q_{\max}$ such that $y_i=1$ and $y_j\leq 1$ for $j=1,\ldots,i-1,i+1,\ldots,n.$ 

We have $\mathcal{D}(\mathcal{H})^{-1}Q(\mathcal{H})\mathcal{D}(\mathcal{H})y=q_{\max}y.$ Comparing $i^{th}$ component both side, we have

$$
d_iy_i+\sum\limits_{j\sim i} \frac{d_j}{d_i}\frac{d_{ij}}{k-1}y_j=q_{\max}y_i.$$ 

Since $y_j\leq y_i=1$, we have 
\begin{equation}\label{eq1}
q_{\max}\leq d_i+\frac{s_i}{k-1}.
\end{equation}

Hence $q_{\max}(\mathcal{H})\leq \max\limits_{i\in V}(d_i+\frac{s_i}{k-1}).$

We now prove the equality part. Suppose $q_{\max}(\mathcal{H})= \max\limits_{i\in V}(d_i+\frac{s_i}{k-1}).$ This implies $q_{\max}\geq d_i+\frac{s_i}{k-1}$ for $i=1,\ldots,n.$ Using Equation \ref{eq1}, we have $q_{\max}=d_i+\frac{s_i}{k-1}.$ This implies $y_j=1$ for $j\sim i.$

Let $j\in V$ and $j\sim i.$ Then we have 

$$
d_jy_j+\sum\limits_{p\sim j} \frac{d_p}{d_j}\frac{d_{pj}}{k-1}y_p=q_{\max}y_j.$$ 

We have noticed that $y_j=1$ for $j\sim i.$ Then $q_{\max}\leq d_j+\frac{s_j}{k-1}.$ from the equality we get $q_{\max}= d_j+\frac{s_j}{k-1}.$ This implies $y_p=1$ for $p\sim j$ and $j\sim i.$ Continuing this process we have $q_{\max}=d_i+\frac{s_i}{k-1}$ for $i=1,\ldots,n.$ Hence $d_1+\frac{s_1}{k-1}=d_2+\frac{s_2}{k-1}=\cdots=d_n+\frac{s_n}{k-1}.$

\textbf{Claim} Let $\mathcal{H}$ be a connected $k$-uniform hypergraph. Then $d_1+\frac{s_1}{k-1}=d_2+\frac{s_2}{k-1}=\cdots=d_n+\frac{s_n}{k-1}$ if and only if either $\mathcal{H}$ is a regular $k$-uniform hypergraph or a bipartite semiregular graph.

\textbf{Proof:} We first assume that $d_1+\frac{s_1}{k-1}=d_2+\frac{s_2}{k-1}=\cdots=d_n+\frac{s_n}{k-1}.$ Let $\delta=d_{n}$ be the minimum degree of $\mathcal{H}$ and let $n$ be a vertex with lowest degree. Let $i_1,\ldots,i_t$ be the vertex adjacent with $n.$ We now show that $i_1,\ldots,i_t$ are of equal degrees. Suppose not, then let $d$ be the highest degree corresponding to the vertex $i_1.$ For vertex $n,$ $\delta+\frac{s_n}{k-1}=\delta+\sum\limits_{k\sim n} \frac{d_k}{d_n}\frac{d_{kn}}{k-1}<\delta+\frac{d}{d_n}\sum\limits_{k\sim n}\frac{d_{kn}}{k-1}=\delta+d.$ For vertex $i_1$, $d+\frac{s_{i_1}}{k-1}\geq d+\delta.$ Therefore $\delta+\frac{s_n}{k-1}<d+\frac{s_{i_1}}{k-1}.$ But this is not possible. Hence $i_1,\ldots,i_t$ have equal degree $d$ (say). Let $i_1$ be adjacent to $n,j_1,\ldots,j_l.$ Since $\delta$ is the lowest degree, then the degree of the vertices $n,j_1,\ldots,j_l$ are greater or equal to $\delta .$ If possible one of them be greater than $\delta.$ Then for vertex $j_1$, $d_{j_1}+\frac{s_{j_1}}{k-1}>d+\delta.$ Therefore $d_n+\frac{s_n}{k-1}=\delta+d<d_{j_1}+\frac{s_{j_1}}{k-1}.$ This is not possible. Hence the degrees of $j_1,\ldots,j_l$ are equal to $\delta.$ In this way we can show that $\delta$ degree vertices are adjacent to $d$ degree vertices and $d$ degree vertices are adjacent to $\delta$ degree vertices.

We now consider the the underlying weighted graph $G$ such that $Q(\mathcal{H})$ is the signless Laplacian matrix of $G$. Suppose $G$ is not bipartite. Then there is an odd cycle, say, $jj_1j_2\cdots j_{k-1}j_{k}j_{k}'j_{k-1}'\cdots j'_2j'_1j.$ Then the degree of the vertex $j$ is either $d$ or $\delta.$ We assume that the degree of the vertex $j$ is $d$. Since $j_1$ and $j'_1$ are adjacent to $j$, then the degrees of $j_1$ and $j'_1$ are $\delta.$ Similarly we can show that the degree of $j_k$ and $j'_k$ are $d$ if $k$ is even and $\delta$ if $k$ is odd. But $j_k$ and $j'_k$ are adjacent vertices. Hence $d=\delta.$ This implies $\mathcal{H}$ is regular. 

If $G$ is bipartite, then $Q(\mathcal{H})$ can be written in the following form $\begin{bmatrix} O & B^t\\ B & 0\end{bmatrix}.$ This is possible only when $k=2.$ Hence in this case $\mathcal{H}$ is a bipartite semiregular graph.\qe

We proved that if equality holds then either $\mathcal{H}$ is a regular hypergraph or a bipartite semiregular graph. 

Converse is very simple.
\qe


\begin{lemma} {\rm\cite{ellingham2000}}
Let $M$ be a real symmetric $n\times n$ matrix, and let $\lambda$ be an eigenvalue of $M$ with an eigenvector $x$ all of whose entries are nonnegative. Let $R_{i}(M)$ be the $i$th row sum of $M$. Then 
\begin{center}
    $\min\limits_{1\leq i\leq n}R_{i}(M)\leq \lambda\leq\max\limits_{1\leq i\leq n}R_{i}(M)$.
\end{center}
Moreover, if the rowsums of M are not all equal and if all entries of x are positive, then both inequalities above are strict.

\end{lemma}

The following theorem is an immediate application of above lemma. 

\begin{theorem}\label{thm2.12}
Let $\mathcal{H}$ be a hypergraph and $d_{\max}$ be the maximum degree of $\mathcal{H}$. Then $2d_{\min}\leq q_{\max}(\mathcal{H})\leq 2d_{\max}$. Equality holds if and only if $\mathcal{H}$ is a regular hypergraph.
\end{theorem}

\begin{lemma} \label{newlem1}
Let $\cal{H}$ be a hypergraph with a signless Laplacian matrix $Q(\mathcal{H})$. Then
\begin{center}
    $\min\limits_{u\in V}(2d^2_u+\frac{2}{k-1}\sum\limits_{v\sim u}d_{uv}d_v)^{\frac{1}{2}}\leq q_{\max}(\mathcal{H})\leq\max\limits_{u\in V}(2d^2_u+\frac{2}{k-1}\sum\limits_{v\sim u}d_{uv}d_v)^{\frac{1}{2}}.$
\end{center}
\end{lemma}

\pr Since $Q({\mathcal{H}})=\mathcal{D}({\mathcal{H}})+\mathcal{A}({\mathcal{H}})$, we have
$R_{u}(Q({\mathcal{H}}))=2d_{u}$ and  $R_{u}(\mathcal{A}^{2}({\mathcal{H}}))=R_{u}(\mathcal{A}({\mathcal{H}})\mathcal{D}({\mathcal{H}}))=\frac{1}{(k-1)}\sum\limits_{u\sim v}d_{uv}d_{v}.$ Now for any $u\in V$ we have,

\begin{align*}
R_{u}(Q^{2}(\mathcal{H}))&=R_{u}\bigg(\big(\mathcal{D}({\mathcal{H}})+\mathcal{A}({\mathcal{H}})\big)\big(\mathcal{D}({\mathcal{H}})+\mathcal{A}({\mathcal{H}})\big)\bigg)\\
&=R_{u}\bigg(\mathcal{D}({\mathcal{H}})\big(\mathcal{D}({\mathcal{H}})+\mathcal{A}({\mathcal{H}})\big)\bigg)+R_{u}\big(\mathcal{A}({\mathcal{H}})\mathcal{D}({\mathcal{H}})\big)+R_{u}\big(\mathcal{A}^{2}({\mathcal{H}})\big)\\
&=d_{u}R_{u}(Q({\mathcal{H}}))+2R_{u}(\mathcal{A}({\mathcal{H}})\mathcal{D}({\mathcal{H}}))\\
&=d_{u}R_{u}\big(Q({\mathcal{H}})\big)+2\frac{1}{(k-1)}\sum\limits_{u\sim v}d_{uv}d_{v}\\
&=2d^2_{u}+2\frac{1}{(k-1)}\sum\limits_{u\sim v}d_{uv}d_{v}.
\end{align*}
Consider the polynomial $P(x)=x^2$. Using Lemma \ref{lem3.7}, we have 
\begin{center}
    $\min\limits_{u\in V}(2d^2_u+\frac{2}{k-1}\sum\limits_{v\sim u}d_{uv}d_v)^{\frac{1}{2}}\leq q_{\max}(\mathcal{H}))\leq\max\limits_{u\in V}(2d^2_u+\frac{2}{k-1}\sum\limits_{v\sim u}d_{uv}d_v)^{\frac{1}{2}}.$
\end{center}\qe
\begin{theorem}
Let $\mathcal{H}$ be a connected $k$-uniform hypergraph with $n(\geq 2)$ vertices and $m$ hyperedges. Then
\begin{center}
    $q_{\max}(\mathcal{H})\leq \frac{\bigg((k-1)(2d^2_{\max}+2md^2_{\min})+2m(km-nd_{\min})\bigg)^{\frac{1}{2}}}{\sqrt{k-1}}$,
\end{center}
where $d_{\max}$ and $d_{\min}$ are the maximum and minimun degrees of $\mathcal{H}$, respectively.
\end{theorem}

\pr From Lemma \ref{newlem1}, we have

\begin{align*}
q_{\max}(\mathcal{H})&\leq\max\limits_{u\in V}(2d^2_u+\frac{2}{k-1}\sum\limits_{v\sim u}d_{uv}d_v)^{\frac{1}{2}}\\
&\leq \max\limits_{u\in V}\bigg(2d^2_{\max}+\frac{2m}{k-1}\big(\sum\limits_{i\in V}d_{i}-d_{u}-\sum\limits_{u(\neq v)\nsim v}d_{v}\big)\bigg)^{\frac{1}{2}}\\
&\leq \max\limits_{u\in V}\bigg(2d^2_{\max}+\frac{2km^{2}}{k-1}-\frac{2md_u}{k-1}-\frac{2md_{\min}(n-1)}{k-1}+2md_{\min} d_u\bigg)^{\frac{1}{2}}\\
&\leq\bigg(2d^2_{\max}+\frac{2km^{2}}{k-1}-\frac{2md_{\min}}{k-1}-\frac{2md_{\min}(n-1)}{k-1}+2md^2_{\min}\bigg)^{\frac{1}{2}}\\
&\leq\frac{\bigg((k-1)(2d^2_{\max}+2md^2_{\min})-2md_{\min}+2m(km-d_{\min}(n-1))\bigg)^{\frac{1}{2}}}{\sqrt{k-1}}\\
&=\frac{\bigg((k-1)(2d^2_{\max}+2md^2_{\min})+2m(km-nd_{\min})\bigg)^{\frac{1}{2}}}{\sqrt{k-1}}.
\end{align*}

Hence the proof.\qe

\begin{theorem}
Let $\mathcal{H}$ be a connected $k$-uniform hypergraph with $n(\geq 2)$ vertices and $m$ hyperedges. If $d_{\max}\leq\frac{km}{n-1}$ then
\begin{center}
    $q_{\max}(\mathcal{H})\geq\frac{\bigg(2(k-1)d^2_{\min}+2km-2d_{\max}(n+2-k-d_{\min})\bigg)^{\frac{1}{2}}}{\sqrt{k-1}},$
\end{center}
where $d_{\max}$ and $d_{\min}$ are the maximum and minimun degrees of $\mathcal{H}$, respectively.
\end{theorem}

\pr From Lemma \ref{newlem1}, we have

\begin{align*}
\hspace{1cm}q_{\max}(\mathcal{H})&\geq\min\limits_{u\in V}(2d^2_u+\frac{2}{k-1}\sum\limits_{v\sim u}d_{uv}d_v)^{\frac{1}{2}}\\
&\geq \min\limits_{u\in V}\bigg(2d^2_{\min}+\frac{2}{k-1}\big(\sum\limits_{i\in V}d_{i}-d_{u}-\sum\limits_{u(\neq v)\nsim v}d_{v}\big)\bigg)^{\frac{1}{2}}\\
&= \min\limits_{u\in V}\bigg(2d^2_{\min}+\frac{2km}{k-1}-\frac{2d_u}{k-1}-\frac{2}{k-1}\sum\limits_{u(\neq v)\nsim v}d_{v}\bigg)^{\frac{1}{2}}\\
&\geq\bigg(2d^2_{\min}+\frac{2km}{k-1}-\frac{2d_{\max}}{k-1}(n+2-k-d_{\min})\bigg)^{\frac{1}{2}},\big(\mbox{ provided }, (n-1)d_{\max}\leq km\big)\\ 
&=\frac{\bigg(2(k-1)d^2_{\min}+2km-2d_{\max}(n+2-k-d_{\min})\bigg)^{\frac{1}{2}}}{\sqrt{k-1}}, \big(\mbox{ provided }, (n-1)d_{\max}\leq km\big).\\
\end{align*}

Hence the proof.\qe

\begin{remark}
    The term $2(k-1)d^2_{\min}+2km-2d_{\max}(n+2-k-d_{\min})$ is not always positive for an arbitrary $k$-uniform hypergraph with $n$ vertices and $m$ hyperedges. Consider the hypergraph $\mathcal{H}=(V,E)$ where $V=\{1,2,3,4,5,6,7\}$ and $E=\big\{\{1,2,3\},\{3,4,5\},\{3,6,7\} \big\}$. Then $2(k-1)d^2_{\min}+2km-2d_{\max}(n+2-k-d_{\min})<0$. For the class of hypergraph $\mathcal{H}$ which satisfies the condition $(n-1)d_{\max}\leq km$, the term $2(k-1)d^2_{\min}+2km-2d_{\max}(n+2-k-d_{\min})$ is always positive. There are many $k$-uniform hypergraphs satisfying the above condition. The following is an example of such hypergraph. Consider $\mathcal{H}=(V,E)$ where $V=\{1,2,3,4,5\}$ and $E=\bigg\{\{1,2,3\},\{1,4,5\},\{3,4,5\}\bigg\}$ are satisfying the condition $(n-1)d_{\max}\leq km$.
    \end{remark}

\begin{remark}
Let $\mathcal{H}=(V,E)$ be a hypergraph on $n$ vertices. Let $x=(x_1,\ldots,x_n)^t$ be a vercor in $\mathbb{R}^n$.
Then $(\mathcal{L}(\mathcal{H})x)(i)=\sum\limits_{i\sim j}\sum\limits_{e\in E_{ij}}\frac{1}{|e|-1}(x(i)-x(j))=d_{i}(\mathcal{H})x(i)-\sum\limits_{j,i\sim j}\sum\limits_{e\in E_{ij}}\frac{1}{|e|-1}x(j)$
\end{remark}

\begin{theorem}
Let $\mathcal{H}=(V,E)$ be a noncomplete connected $k$-uniform hypergraph with $n$ vertex. Let $X$ and $Y$ be two subsets of $V$ such that $X\cap Y=\phi$ and the vertices of $X$ are not adjacent with the vertices of $Y.$ Then
\begin{center}
   $q_{\max}(\mathcal{H})\geq\frac{\beta^2_{1}|X|+\beta^2_{2}|Y|}{k-1}$,
\end{center}
where $\beta_{1}$ and $\beta_{2}$ are two real numbers such that $\beta^2_{1}+\beta^2_{2}=1$.
\end{theorem}

\pr Let $V$ be the vertex set of $\mathcal{H}$. Then $V=X\cup Y\cup V\backslash(X\cup Y)$. 

Let us consider the vector $x=(x_1,\ldots,x_n)^t$ in $\mathbb{R}^n$ such that,

$$
x(i)=\begin{cases}
    \beta_{1}~~~,\mbox{ if }~~i\in X,\\
     -\beta_{2}~~~,\mbox{ if }~~i\in Y,\\
     0~~~~~~,\mbox{ elsewhere }
 \end{cases}$$

Let $n_{i}(\mathcal{H})=\sum\limits_{\underset{i\sim j}{j\in V\backslash(X\cup Y)}}\sum\limits_{e\in E_{ij}}\frac{1}{|e|-1}$, for all $i\in V\backslash(X\cup Y)$.

Now $\forall i\in X$,
\begin{align*}
(\mathcal{L}(\mathcal{H})x)(i)&=d_{i}(\mathcal{H})x(i)-\sum\limits_{\underset{j\sim i}{j\in X}}\sum\limits_{e\in E_{ij}}\frac{1}{|e|-1}x(j)\\
&=\beta_{1} d_{i}(\mathcal{H})-\beta_{1}\sum\limits_{\underset{j\sim i}{j\in X}}\sum\limits_{e\in E_{ij}}\frac{1}{|e|-1}.\\
&=\beta_{1} d_{i}(\mathcal{H})-\beta_{1}\big(d_{i}(\mathcal{H})-n_{i}(\mathcal{H})\big).
\end{align*}

$\therefore (\mathcal{L}(\mathcal{H})x)(i)=\beta_{1} n_{i}(\mathcal{H})$, for all $i\in X$.
Similarly, $(\mathcal{L}(\mathcal{H})x)(i)=-\beta_{2} n_{i}(\mathcal{H})$, for all $i\in Y$.

Now 
\begin{align*}
\mu_{max}(\mathcal{H})&\geq\langle(L(\mathcal{H})x,x\rangle\\
&=\sum\limits_{i\in V}(L(\mathcal{H})x)(i)x(i)\\
&=\sum\limits_{s\in X}\beta^2_{1}n_{i}(\mathcal{H})+\sum\limits_{s\in Y}\beta^2_{2}n_{i}(\mathcal{H})\\
&\geq\frac{\beta^2_{1}|X|+\beta^2_{2}|Y|}{k-1}.
\end{align*}

Hence using Theorem $\ref{lem2.7}$ we get, $q_{\max}(\mathcal{H})\geq\frac{\beta^2_{1}|X|+\beta^2_{2}|Y|}{k-1}$.\qe

\section{Smallest signless Laplacian eigenvalue of a $k$-uniform hypergraph}
The following result was first proved in 1994 by Desai and Rao \cite{Desai1994}. 
\begin{lemma}\label{thm:1}
Let $G$ be a connected, positive weighted graph. Then the smallest signless Laplacian eigenvalue is $0$ if and only if $G$ is bipartite. 
\end{lemma}

It is natural to ask to characterize the $k$-uniform hypergraphs such that the smallest signless Laplacian eigenvalue is $0$. The following theorem answers this question.

\begin{theorem}\label{t2.13}
Let $\mathcal{H}$ be a connected $k(\geq 3)$-uniform hypergraph. Then the smallest signless Laplacian eigenvalue of $\mathcal{H}$ cannot be zero. \end{theorem}

\pr Suppose that the smallest signless Laplacian eigenvalue of $\mathcal{H}$ is zero. Let $x$ be an eigenvector corresponding to zero.
Let $G$ be a weighted graph such that the matrix $Q(\mathcal{H})$ is the signless Laplacian matrix of $G$. Using Lemma \ref{thm:1}, $G$ is a bipartite graph. Let $V_1$ and $V_2$ be two disjoint partitions of the vertex set $V$. Let $i$ and $j$ be two vertices in $G$ such that $i\sim j$. Therefore $i\in V_1$ and $j\in V_2$. Since $i\sim j$ in $G$, we have $i\sim j$ in $\mathcal{H}$. Let $\{i,i_1,\ldots,i_{k-2}, j\}$ be a hyperedge in $\mathcal{H}$, since $k\geq3$. This says that $i\sim i_1$ and $i_1\sim j$ in $\mathcal{H}$.
Therefore $i\sim i_1$ and $i_1\sim j$ in $G$ which is a contradiction as $G$ is bipartite. Hence the smallest signless Laplacian eigenvalue of $\mathcal{H}$ cannot not be zero.\qe
 
\begin{theorem}
Let $\mathcal{H}$ be a connected $k$-uniform hypergraph. Then following are true.
\begin{enumerate}
\item If $k=2$, then the smallest positive eigenvalue of $\mathcal{H}$ is zero if and only if $\mathcal{H}$ is bipartite.

\item If $k\geq 3$, then the smallest positive eigenvalue of $\mathcal{H}$ cannot be zero. 
\end{enumerate}
\end{theorem}

\begin{theorem}\label{t2.12}
Let $\mathcal{H}$ be a connected $k$-uniform hypergraph with at least two vertices. Let $u$ and $v$ be two distinct vertices. Then $q_{\min}(\mathcal{H})\leq\frac{d_u+d_v}{2}$
\end{theorem}

\pr Let V be the vertex set and $u,v\in V$. 
Consider the vector $x=(x_1,\ldots,x_n)^t$ such that 

\begin{center}
$x(i) = \left\{ \begin{array}{rcl}
\frac{1}{\sqrt{2}} & \mbox{if}
& i\sim u \\ -\frac{1}{\sqrt{2}} & \mbox{if} & i\sim v \\
0 & else
\end{array}\right.$
\end{center}

Now  $q_{\min}(\mathcal{H})=\min\limits_{||x||=1} x^t Q(\mathcal{H})x$

Therefore 
\begin{align*}
q_{\min}(\mathcal{H})&\leq x^t Q(\mathcal{H})x\\ &=\sum\limits_{i\sim j}q_{ij}(x(i)+x(j))^2\\
&\leq\sum\limits_{i\sim u}q_{ij}(\frac{1}{\sqrt{2}})^2+\sum\limits_{j\sim v}q_{ij}(-\frac{1}{\sqrt{2}})^2\\
&\leq \frac{d_u+d_v}{2}.
\end{align*}
\qe

\begin{theorem}
Let $\mathcal{H}=(V,E)$ be a connected $k$-uniform hypergraph with at least two vertices and let $X$ and $Y$ be two subsets of $V$ such that $X\cap Y=\phi.$ Assume that the vertices of $X$ are not adjacent with the vertices of $Y.$ Then
\begin{center}
   $q_{\min}(\mathcal{H})\leq\frac{(|X|+|Y|)m}{2(k-1)}$.
\end{center}
\end{theorem}

\pr Consider the vector $x=(x_1,\ldots,x_n)^t$ such that

$$
x(i)==\begin{cases}
    \frac{1}{\sqrt{2}}~~~,\mbox{ if }~~i\in X,\\
     -\frac{1}{\sqrt{2}}~~~,\mbox{ if }~~~i\in Y,\\
     0~~~~~~~~,\mbox{ elsewhere }
 \end{cases}$$

Now 
\begin{align*}
q_{\min}(\mathcal{H})&\leq x^tQ(\mathcal{H})x\\
&=\sum\limits_{i\sim j}\sum\limits_{e\in E_{ij}}\frac{1}{|e|-1}(x(i)+x(j))^2\\
&=\sum\limits_{\underset{i\sim s}{i\in X}}\sum\limits_{e\in E_{is}}\frac{1}{|e|-1}x^2(i)+\sum\limits_{\underset{j\sim t}{j\in Y}}\sum\limits_{e\in E_{jt}}\frac{1}{|e|-1}x^2(j)\\
&\leq\frac{(|X|+|Y|)m}{2(k-1)}
\end{align*}
Hence it follows.\qe

Using Theorem \ref{thm:2.6} and Lemma \ref{lem2.7}, we have $q_{\max}\geq d_{\max}+\frac{1}{k-1}.$ The following theorem supplies a similar type of upper bounds for $q_{\min}.$ 

\begin{theorem}\label{cor1}
Let $\mathcal{H}$ be a connected $k$-uniform hypergraph with at least two vertices. 
Then \begin{center}
    $q_{\min}(\mathcal{H})\leq d_{\max}-\frac{1}{k-1}.$
\end{center}

\end{theorem}

\pr Let $i$ and $j$ be two adjacent vertices in $\mathcal{H}$ and let $d_i$ and $d_j$ be the degrees of the vertices $i$ and $j$, respectively.

Let $x=(x_1,\ldots,x_n)^t$ be any vector in $\mathbb{R}^n$. Then $x^tQ(\mathcal{H})x=\sum\limits_{i\sim j}\sum\limits_{e\in E_{ij}}\frac{1}{|e|-1}(x(i)+x(j))^2.$

We have $q_{\min}(\mathcal{H})\leq\frac{\sum\limits_{i\sim j}\sum\limits_{e\in E_{ij}}\frac{1}{|e|-1}(x(i)+x(j))^2}{\sum\limits_{i=1}^n x^2_i}.$

Let $N_i$ and $N_j$ be the cardinality of the set $\{t : t\sim i\}$ and
$\{t : t\sim j\},$ respectively. Let $d_i$ and $d_j$ be the degrees of $i$ and $j$, respectively. 
There are $d_i$ hyperedges containing the vertex $i$ and each set contains atmost $k-1$ vertices which are adjacent to $i$. Then $N_i\leq (k-1)d_i.$ Similarly $N_j\leq (k-1)d_j.$

Let $x=(0,0,\ldots,0,\underbrace{1}_{ith},0,\ldots,0,\underbrace{-1}_{jth},0,\ldots,0).$ Then $\frac{\sum\limits_{i\sim j}\sum\limits_{e\in E_{ij}}\frac{1}{|e|-1}(x(i)+x(j))^2}{\sum\limits_{i=1}^n x^2_i}=\frac{N_i+N_j-2}{2(k-1)}\leq \frac{(k-1)d_i+(k-1)d_j-2}{2(k-1)}.$ Hence $q_{\min}(\mathcal{H})\leq \frac{d_i+d_j-\frac{2}{k-1}}{2}.$
This implies $q_{\min}(\mathcal{H})\leq d_{\max}-\frac{1}{k-1}.$
\qe

\begin{theorem}
Let $\mathcal{H}$ be a connected $k$-uniform hypergraph with $n(\geq 2)$ vertices and $m$ hyperedges. 
Then \begin{center}
    $q_{\min}(\mathcal{H})\leq\sqrt{2d^2_{\max}+\frac{2km}{k-1}\bigg(km-(n-1)d_{\min}+((k-1)d_{\min}-1)d_{\max}\bigg)}$,
\end{center}

where $d_{\max}$ and $d_{\min}$ are two maximum and minimum degees of $\mathcal{H}$, respectively.
\end{theorem}

\pr Let $x=(x_1,\ldots,x_n)^t$ be an unit eigenvector corresponding to $q_{\min}(\mathcal{H})$. 

Then for $i\in V$, $q_{\min}(\mathcal{H})x_{i}=\sum\limits_{j\sim i}q_{ij}(x_i+x_j)$.

Now applying the Cauchy–Schwarz inequality we get,

$\hspace{3cm}q_{\min}(\mathcal{H})^2x^2_{i}\leq\sum\limits_{j\sim i}q_{ij}\sum\limits_{j\sim i}q_{ij}(x_i+x_j)^2=d_i\sum\limits_{j\sim i}q_{ij}(x_i+x_j)^2$

$\hspace{2.5cm}\Rightarrow q_{\min}(\mathcal{H})^2x^2_{i}\leq2d_i\sum\limits_{j\sim i}q_{ij}(x^2_i+x^2_j)$

$\hspace{5cm}\leq2d^2_ix^2_i+2d_i\sum\limits_{j\sim i}q_{ij}x^2_j$

Now taking sum over all the vertices we get, 

$\hspace{2cm}q_{\min}(\mathcal{H})^2\sum\limits_{i\in V}x^2_{i}\leq2\sum\limits_{i\in V}d^2_{i}x^2_{i}+2\sum\limits_{i\in V}d_i(\sum\limits_{j\sim i}q_{ij}x^2_j)$

$\hspace{2.5cm}\Rightarrow q_{\min}(\mathcal{H})^2\leq2d^2_{\max}+2\sum\limits_{i\in V}d_i(\sum\limits_{j\sim i}q_{ij}x^2_j)$ 

$\hspace{4.5cm}\leq2d^2_{\max}+\frac{2m}{k-1}\sum\limits_{i\in V}d_i(\sum\limits_{j\sim i}x^2_j)$ 

$\hspace{4.5cm}=2d^2_{\max}+\frac{2m}{k-1}\sum\limits_{i\in V}d_i\big(1-x^2_i-\sum\limits_{\underset{j\neq i}{j\nsim i}}x^2_j\big)$ 

$\hspace{4.5cm}=2d^2_{\max}+\frac{2km^2}{k-1}-\frac{2m}{k-1}\bigg(\sum\limits_{i\in V}d_ix^2_i+\sum\limits_{i\in V}d_i\sum\limits_{\underset{j\neq i}{j\nsim i}}x^2_j\bigg)$

$\hspace{4.5cm}\leq2d^2_{\max}+\frac{2km^2}{k-1}-\frac{2m}{k-1}\bigg(\sum\limits_{i\in V}d_ix^2_i+\sum\limits_{i\in V}d_{\min}(n-1-(k-1)d_i)x^2_i\bigg)$

$\hspace{4.5cm}=2d^2_{\max}+\frac{2km^2}{k-1}-\frac{2m}{k-1}\bigg((n-1)d_{\min}-((k-1)d_{\min}-1)\sum\limits_{i\in V}d_ix^2_i\bigg)$

Hence $q_{\min}(\mathcal{H})\leq\sqrt{2d^2_{\max}+\frac{2km}{k-1}\bigg(km-(n-1)d_{\min}+((k-1)d_{\min}-1)d_{\max}\bigg)}$.

The following theorem gives an upper bound for $q_{\min}$ in terms of first Zagreb index.
\begin{theorem}
Let $\mathcal{H}$ be a connected hypergraph of order $n(\geq 2)$. Then \begin{center}
    $q_{\min}(\mathcal{H})\leq2\sqrt{\frac{Z_{1}(\mathcal{H})}{n}}$,
\end{center}
where $Z_1(\mathcal{H})$ is first Zagreb index.
\end{theorem}

\pr Let $x=(x_1,x_2,\ldots,x_n)$ be a unit vector in $\mathbb{R}^n.$ Let $\beta=\max\limits_{1\leq i\leq n} |x_i|.$
Since for $i\in V$, $q_{\min}(\mathcal{H})x_{i}=\sum\limits_{j\sim i}q_{ij}(x_i+x_j)$ and applying the Cauchy–Schwarz inequality we get, 

$\hspace{3cm}q_{\min}(\mathcal{H})^2x^2_{i}\leq\sum\limits_{j\sim i}q_{ij}\sum\limits_{j\sim i}q_{ij}(x_i+x_j)^2=d_i\sum\limits_{j\sim i}q_{ij}(x_i+x_j)^2$

$\hspace{2.5cm}\Rightarrow q_{\min}(\mathcal{H})^2x^2_{i}\leq d_i\sum\limits_{i\sim j}q_{ij}(x_i+x_j)^2\leq4\beta^2 d^2_{i}$. 
 
 After taking sum over all the vertices we get, $q_{\min}(\mathcal{H})\leq2\beta\sqrt{Z_{1}(\mathcal{H})}$. Now take $x=(\frac{1}{\sqrt{n}},\frac{1}{\sqrt{n}},\ldots,\frac{1}{\sqrt{n}}).$ Here $\beta=\frac{1}{\sqrt{n}}.$ Hence $q_{\min}(\mathcal{H})\leq2\sqrt{\frac{Z_{1}(\mathcal{H})}{n}}.$
\qe

\section{Signless Laplacian spread of a hypergraph}
For any $n\times n$ complex matrix $M$, the spread $s(M)$ of $M$ is defined as the diameter of its spectrum, that is, $s(M)=\max\limits_{i,j}|\lambda_i-\lambda_j|$, where the maximum is taken over all pair of eigenvalues of $M.$ There is a considerable literature on the spread of an arbitrary
matrix \cite{johson1985}. Let $M$ be a real symmetric matrix and let $\lambda_{\max}(M)=\lambda_1(M)\geq\lambda_2(M)\geq\cdots\geq \lambda_n(M)=\lambda_{\min}(M)$ be the eigenvalues of $M.$ Then the spread of $M$ is defined by 
\begin{center}
    $s(M)=\lambda_{\max}(M)-\lambda_{\min}(M)$
\end{center}

There are several results concerning the spread of the adjacency matrix, Laplacian matrix and signless Laplacian matrix of a graph, see for example \cite{fan2008, johson1985, liu2010} and references therein.

Motivated from the above definition one can define the signless Laplacian spread of hypergraphs. 
The signless Laplacian spread of a hypergraphs is defined as 
\begin{center}
    $s_{Q}(\mathcal{H})=q_{\max}(\mathcal{\mathcal{H}})-q_{\min}(\mathcal{H})$,
\end{center}
where $q_{\max}(\mathcal{H})$ and $q_{\min}(\mathcal{H})$ are the maximum and minimum signless Laplacian eigenvalues of $\mathcal{H},$ respectively.

In this section, we obtain bounds for the signless Laplacian spread $s_Q(\mathcal{H})$ of $k$-uniform hypergraphs.

\subsection{Signless Laplacian spread of a hypergraphs in terms of maximum, minimum degrees}

In  \cite{Oliveira2010}, the authors proved that the signless Laplacian spread $s_Q(G)\geq 2$ for any graph. But for hypergraphs this is not a lower bound. For example, consider the $3$-uniform  hypergraph $\mathcal{H}$ on three vertices. We can easily check that the signless Laplacian spread of $\mathcal{H}$ is $1.5$ which is smaller than $2.$  
The following theorem says that $Q$-spread of a $k$-uniform hypergraph cannot be less than $1.$

\begin{theorem} \label{t3.3}
Let $\mathcal{H}$ be a hypergraph with $n$ vertices and at least one hyperedge. Then $s_{Q}(\mathcal{H})>1.$ 
\end{theorem}

\pr We first prove that $q_{\min}(\mathcal{H})\leq d_{\min}.$ Let $V(\mathcal{H})=\{1,\ldots,n\}$ and let the degree of the vertex $n$ is $d_{\min}.$ Taking $l=(0,0,,\ldots,1)^t$ is an unit vector which is not an eigenvector corresponding to $d_{\min}$. Now 
\begin{align*}
q_{\min}(\mathcal{H})&\leq l^t Q(\mathcal{H})l\\
&=\sum\limits_{i\sim j}q_{ij}(l_i+l_j)^2\\
&=\sum\limits_{i\sim n}q_{in}=d_{\min}.
\end{align*}

Using Theorem \ref{thm2.12}, we have $q_{\max}\geq 2d_{\min}.$
Therefore $s_{Q}(\mathcal{H})\geq d_{\min}\geq1.$ 

We now show that the inequality is strict. Suppose $s_{Q}(\mathcal{H})=1$ for some $k$-uniform hypergraphs $\mathcal{H}$, then from the proof, we have $q_{\max}=2d_{\min}$ and $q_{\min}=d_{\min}$. Using Theorem $\ref{thm2.12}$, we have $\mathcal{H}$ is a regular $k$-uniform hypergraph. Also $q_{\min}=d_{\min}$ implies $(0,0,0,\ldots,1)^t$ is an eigenvector of $Q$ corresponding the eigenvalue $q_{\min}$. 
If $(0,0,0,\ldots,1)^t$ is an eigenvector, then the $n^{th}$ vertex in $\mathcal{H}$ is isolated. Then the degree of each vertex is zero as $\mathcal{H}$ is regular. A contradiction to the fact that $\mathcal{H}$ has at least one hyperedge. Hence $s_{Q}(\mathcal{H})>1.$\qe 

\begin{proposition}
Let $\mathcal{F}$ be the set of all hypergraphs with finite vertices and at least one hyperedge. Then $\inf\big\{s_Q(\mathcal{H}): \mathcal{H}\in \mathcal{F}\big\}=1.$
\end{proposition}

\pr Using Theorem \ref{t3.3}, a lower bound of the set $\big\{s_Q(\mathcal{H}): \mathcal{H}\in \mathcal{F}\big\}$ is $1$.

We have to prove that the infimum of the set $\big\{s_Q(\mathcal{H}): \mathcal{H}\in\mathcal{F}\big\}$ is $1$. Suppose that the infimum is not $1.$ Let $1+M$ is the infimum where $M$ is a positive real number. Using Archemedian property of real numbers, there exists a natural number $n_1$ such that $\frac{1}{n_1}< M.$ Consider the hypergraph $\mathcal{H}=(V,E)$, where $V=\{1,2,\ldots,n_1+1\}$ and $E=\big\{\{1,2,\ldots,n_1+1\}\big\}$.

Then $Q(\mathcal{H})=\begin{bmatrix}
1 & \frac{1}{n_1} & \frac{1}{n_1} & \cdots & \frac{1}{n_1}\\
\frac{1}{n_1} & 1 & \frac{1}{n_1} & \cdots & \frac{1}{n_1}\\
\vdots & \vdots & \vdots & \vdots & \vdots\\
\frac{1}{n_1} & \frac{1}{n_1} & \frac{1}{n_1} & \cdots & 1
\end{bmatrix}.$

The eigenvalues of $Q(\mathcal{H})$ are $2$ with multiplicity $1$ and $1-\frac{1}{n_1}$ with multiplicity $n_1.$ Then $s_Q(\mathcal{H})=2-1+\frac{1}{n_1-1}=1+\frac{1}{n_1}$ which is less than $1+M.$ A contradiction to the fact that $1+M$ is the infimum of $\big\{s_Q(\mathcal{H}): \mathcal{H}\in \mathcal{F}\big\}.$ Hence the infimum of $\big\{s_Q(\mathcal{H}): \mathcal{H}\in \mathcal{F}\big\}$ is $1$.\qe

The following theorem gives a lower bound for the spread in terms of minimum degree and total number edges. We use the notation $\overline{d}=\frac{d_1+\ldots+d_n}{n}$.

\begin{theorem}\label{thm:2.8}
 Let $\mathcal{H}$ be a $k$-uniform hypergraph with $n$ vertices and $m(\geq 1)$ hyperedges. Assume that $d_{\min}> \frac{\overline{d}}{2}$. Then
 \begin{center}
    $s_Q(\mathcal{H})\geq\frac{2nd_{\min}-km}{n-1}$,
 \end{center}
 Equality holds if and only if the form signless Laplacian matrix is $Q(\mathcal{H})= t((n-1)I_n+J_n)$, where $J_n$ is the $n\times n$ matrix with each entry $1$, $I_n$ is the identity matrix of order $n$ and $t$ is a positive real number.
 \end{theorem}

\pr Since $Tr(Q(\mathcal{H})=\sum\limits_{i=1}^{n} q_{i}(\mathcal{H}).$ Using Theorem $\ref{thm2.5}$, we have  
\begin{align*}
km&=\sum\limits_{i\in V}d_{i}=\sum\limits_{i=1}^{n} q_{i}(\mathcal{H})\\
&\geq q_{\max}(\mathcal{H})+(n-1)q_{\min}(\mathcal{H})
\end{align*}

Therefore $q_{\min}(\mathcal{H})\leq\frac{km-q_{\max}(\mathcal{H})}{n-1}$. Now
\begin{align*}
s_{Q}(\mathcal{H})&=q_{\max}(\mathcal{H})-q_{\min}(\mathcal{H})\\
&\geq q_{\max}(\mathcal{H})-\frac{km-q_{\max}(\mathcal{H})}{n-1}\\
\end{align*}
\begin{equation}\label{neweq1}
\geq\frac{n q_{\max}(\mathcal{H})-km}{n-1}.
\end{equation}
 Using Theorem \ref{thm2.12}, we have $s_Q(\mathcal{H})\geq\frac{2nd_{\min}-km}{n-1}$. As we know that spread is always positive for a $k$-uniform hypergraph with at least one hyperedge. We notice that $\frac{2nd_{\min}-km}{n-1}$ is positive if and only if $d_{\min}> \frac{\overline{d}}{2}$.

We now prove the equality part. If equality holds, then from the above proof, we have $q_2=q_3=\cdots=q_n$ and $q_{\max}(\mathcal{H})=2d_{\min}$ and using Theorem $\ref{thm2.12}$, $\mathcal{H}$ is a regular hypergraph. This implies $(\frac{1}{\sqrt{n}},\ldots,\frac{1}{\sqrt{n}})^t$ is an eigenvector corresponding to $q_1.$

 Notice that the signless Laplacian matrix of $\mathcal{H}$ is positive semi-definite, then the largest eigenvalue of $Q(\mathcal{H})$ is equal to the spectral radius. Since $Q(\mathcal{H})$ is irreducible and nonnegative, the largest eigenvalue is simple. So it has two distinct eigenvalues $q_1$ with multiplicity $1$ and $q_2$ with multiplicity $n-1.$ The matrix $Q(\mathcal{H})$ is orthogonally diagonalizable. Hence there exists a orthogonal matrix $P$ whose first column is $(\frac{1}{\sqrt{n}},\ldots,\frac{1}{\sqrt{n}})^t$ such that 

\begin{align*}
P^tQ(\mathcal{H})P&=\begin{bmatrix}
q_1 & 0 & 0 & \cdots & 0\\
0  & q_2 & 0 & \cdots & 0\\
\vdots & \vdots & \vdots & \cdots & 0\\
0 & 0 & 0 & \cdots & q_2
\end{bmatrix}\\
&=\begin{bmatrix}
q_2 & 0 & 0 & \cdots & 0\\
0  & q_2 & 0 & \cdots & 0\\
\vdots & \vdots & \vdots & \cdots & 0\\
0 & 0 & 0 & \cdots & q_2
\end{bmatrix}+\begin{bmatrix}
q_1-q_2 & 0 & 0 & \cdots & 0\\
0  & 0 & 0 & \cdots & 0\\
\vdots & \vdots & \vdots & \cdots & 0\\
0 & 0 & 0 & \cdots & 0
\end{bmatrix}
\end{align*}
\begin{align*}
Q(\mathcal{H})&=P\begin{bmatrix}
q_2 & 0 & 0 & \cdots & 0\\
0  & q_2 & 0 & \cdots & 0\\
\vdots & \vdots & \vdots & \cdots & 0\\
0 & 0 & 0 & \cdots & q_2
\end{bmatrix}P^t+P\begin{bmatrix}
q_1-q_2 & 0 & 0 & \cdots & 0\\
0  & 0 & 0 & \cdots & 0\\
\vdots & \vdots & \vdots & \cdots & 0\\
0 & 0 & 0 & \cdots & 0
\end{bmatrix}P^t
\end{align*}
\begin{align*}
&=\begin{bmatrix}
q_2 & 0 & 0 & \cdots & 0\\
0  & q_2 & 0 & \cdots & 0\\
\vdots & \vdots & \vdots & \cdots & 0\\
0 & 0 & 0 & \cdots & q_2
\end{bmatrix}+P\begin{bmatrix}
q_1-q_2 & 0 & 0 & \cdots & 0\\
0  & 0 & 0 & \cdots & 0\\
\vdots & \vdots & \vdots & \cdots & 0\\
0 & 0 & 0 & \cdots & 0
\end{bmatrix}P^t
\end{align*}
\begin{align*}
&=q_2\begin{bmatrix}
1 & 0 & 0 & \cdots & 0\\
0  & 1 & 0 & \cdots & 0\\
\vdots & \vdots & \vdots & \cdots & 0\\
0 & 0 & 0 & \cdots & 1
\end{bmatrix}+P\begin{bmatrix}
q_1-q_2 & 0 & 0 & \cdots & 0\\
0  & 0 & 0 & \cdots & 0\\
\vdots & \vdots & \vdots & \cdots & 0\\
0 & 0 & 0 & \cdots & 0
\end{bmatrix}P^t\\
&=q_2I_n+\frac{q_1-q_2}{n} J_n.\\
\end{align*}

The degree of $Q(\mathcal{H})$ is $\frac{q_1}{2}$ which is equal to $q_2+\frac{q_1-q_2}{n}.$ This implies $q_1=t(2n-2)$ and $q_2=t(n-2)$ for some positive real number $t$. Then $Q(\mathcal{H})=t((n-2)I_n+J_n).$

Converse part is simple.
\qe

\begin{proposition}\label{cor:2.9}
    Let $\mathcal{H}$ be a $k$-uniform hypergraph with $n$ vertices and $m(\geq 1)$ hyperedges. Then
 \begin{center}
     $s_{Q}(\mathcal{H})\geq\frac{n(d_{\max}+\frac{1}{k-1})-km}{n-1}$.
 \end{center}
 where $d_{\max}$ is the maximum degree of $\mathcal{H}$.
\end{proposition}

\pr Using Equation \ref{neweq1}, we have $s_Q(\mathcal{H})\geq\frac{nq_{\max}-km}{n-1}$. Also using Theorem $\ref{thm:2.6}$ and Lemma $\ref{lem2.7}$, $q_{\max}\geq d_{\max}+\frac{1}{k-1}$. Hence the proof.\qe

\begin{theorem}\label{thm:2.13}
Let $\mathcal{H}$ be a $k$-uniform hypergraph. Then
 \begin{center}
    $s_{Q}(\mathcal{H})\geq d_{\max}-d_{\min}+\frac{1}{k-1}$,
 \end{center}
where $d_{\max}$ and $d_{\min}$ be the maximum and minimum degrees of $\mathcal{H}$, respectively. 
\end{theorem}

\pr Let $x=(0,0,\ldots,1)^t$ be the unit vector then $q_{\min}(\mathcal{H})\leq d_{\min}$.

Hence using Lemma $\ref{lem2.7}$ and Theorem $\ref{thm:2.6}$ we get, $ s_{Q}(\mathcal{H})\geq d_{\max}-d_{\min}+\frac{1}{k-1}$.\qe

\begin{theorem}
Let $\mathcal{H}$ be a $r$-regular $k$-uniform hypergraph. Then $s_{\mathcal{A}}(\mathcal{H})=s_Q(\mathcal{H}).$ 
\end{theorem}

\pr Let $\lambda_1\geq\cdots\geq\lambda_n$ be the eigenvalues of $\mathcal{A}(\mathcal{H}).$
The signless Laplacian matrix of $\mathcal{H}$ is 

$$
Q(\mathcal{H})=\mathcal{D}(\mathcal{H})+\mathcal{A}(\mathcal{H})=rI_n+\mathcal{A}(\mathcal{H}).
$$

Therefore $q_i=r+\lambda_i$ , for $i=1,\ldots,n$ are the eigenvalues of $Q(\mathcal{H}).$ Hence $s_{\mathcal{A}}(\mathcal{H})=s_Q(\mathcal{H}).$
\qe

Next we supply a upper bound of a signless Laplacian spread using only maximum degree of a $k$-uniform hypergraphs.

\begin{lemma} {\rm\cite{Sedrakyan2018}}
If $x,y\geq\frac{1}{2}$ then $\big(\frac{x^2-y^2}{2}\big)^2+\frac{x+y}{2}\geq\sqrt{\frac{x^2+y^2}{2}}$.
\end{lemma}

\begin{theorem}
Let $\mathcal{H}$ be a connected $k$-uniform hypergraph. Then  
\begin{center}
    $s_{Q}(\mathcal{H})<\frac{1}{2\sqrt{2}}\big[(4d^2_{\max}-\frac{1}{4})^2+2(2d_{\max}+\frac{1}{2})\big]$,
\end{center}
where $d_{\max}$ is the maximum degree of $\mathcal{H}$.
\end{theorem}

\pr Let $x=s_{Q}(\mathcal{H})$ and $y=q_{\min}(\mathcal{H})+\frac{1}{2}$.
Therefore by Theorem $\ref{t3.3}$, $x\geq\frac{1}{2}$ and Theorem $\ref{t2.13}$, $y\geq\frac{1}{2}$. 
From above inequality we get,
\begin{align*}
\sqrt{\frac{s_{Q}(\mathcal{H})^2+\big(q_{\min}(\mathcal{H})+\frac{1}{2}\big)^2}{2}}&\leq\bigg(\frac{s_{Q}(\mathcal{H})^2-(q_{\min}(\mathcal{H})+\frac{1}{2})^2}{2}\bigg)^2+\frac{s_{Q}(\mathcal{H})+q_{\min}(\mathcal{H})+\frac{1}{2}}{2}\\
\Rightarrow\frac{s_{Q}(\mathcal{\mathcal{H}})}{\sqrt{2}}&<\bigg(\frac{(q_{\max}(\mathcal{H})+\frac{1}{2})(q_{\max}(\mathcal{H})-2q_{\min}(\mathcal{H})-\frac{1}{2}}{2}\bigg)^2+\frac{q_{\max}(\mathcal{H})+\frac{1}{2}}{2}\\
\Rightarrow s_{Q}(\mathcal{\mathcal{H}})&<\sqrt{2}\bigg[\bigg(\frac{(q_{\max}(\mathcal{H})+\frac{1}{2})(q_{\max}(\mathcal{H})-\frac{1}{2})}{2}\bigg)^2+\frac{q_{\max}(\mathcal{H})+\frac{1}{2}}{2}\bigg]
\end{align*}

Therefore $s_{Q}(\mathcal{H})<\sqrt{2}\bigg[\frac{(q_{\max}(\mathcal{H})^2-\frac{1}{4})^2}{4}+\frac{(q_{\max}(\mathcal{H})+\frac{1}{2})}{2}\bigg]$

Hence using Theorem $\ref{thm2.12}$ we get, \begin{center}
   $s_{Q}(\mathcal{H})<\frac{1}{2\sqrt{2}}\big[(4d^2_{\max}-\frac{1}{4})^2+2(2d_{\max}+\frac{1}{2})\big]$
\end{center}. \qe

\begin{corollary}
Let $\mathcal{H}$ be a connected $r$-regular $k$-uniform hypergraph. Then
 \begin{center}
   $s_{Q}(\mathcal{H})<\frac{1}{2\sqrt{2}}\big[(4r^2-\frac{1}{4})^2+2(2r+\frac{1}{2})\big]$.
 \end{center}
\end{corollary}

\subsection{Signless Laplacian spread of hypergraphs in terms of weak independent number}

\begin{lemma}\rm\cite{Dragomir1994}\label{l3.27}
If $a=(a_1,\ldots,a_n)$ and $b=(b_1,\ldots,b_n)$ are sequences of nonnegative real numbers such that $a^p_1-a^p_2-\ldots-a^p_n\geq0$ and $b^p_1-b^p_2-\ldots-b^p_n\geq0$ then for $0\leq p\leq2$, \begin{center}
    $(a^p_1-a^p_2-\ldots-a^p_n)(b^p_1-b^p_2-\ldots-b^p_n)\leq(a_1b_1-a_2b_2-\ldots-a_nb_n)^p$
\end{center}
For $p=2$ equality holds iff a and b are proportional. 
\end{lemma}

The following theorem is an upper bound of a Q-spread in terms of maximum degree and weak independent number of $\mathcal{H}$. 

\begin{theorem}
Let $\mathcal{H}$ be a hypergraph and let $\tau(\geq1)$ be the weak independence number of $\mathcal{H}$. Then \begin{center}
    $s_{Q}(\mathcal{H})<\frac{2\tau^2d_{\max}}{(\tau^2-\frac{1}{n})}$,
\end{center}
\end{theorem}

\pr Let $a_1=q_{\max}$, $a_i=\frac{q_{\min}}{\sqrt{n}}$, $i=2,3,\ldots,n+1$ and $b_1=\tau$, $b_i=\frac{1}{n}$, $i=2,3,\ldots,n+1$.

For $p=2$ and using Lemma $\ref{l3.27}$,

$\hspace{2cm}\big(q^2_{\max}\underbrace{-\frac{q^2_{\min}}{n}-\ldots-\frac{q^2_{\min}}{n}}_{n-terms}\big)\big(\tau^2\underbrace{-\frac{1}{n^2}-\ldots-\frac{1}{n^2}}_{n-terms}\big)\leq\big(q_{\max}\tau\underbrace{-\frac{q_{\min}}{n\sqrt{n}}-\ldots-\frac{q_{\min}}{n\sqrt{n}}}_{n-terms}\big)^2$

$\hspace{1.5cm}\Rightarrow(q^2_{\max}-q^2_{\min})(\tau^2-\frac{1}{n})\leq(q_{\max}\tau-\frac{q_{\min}}{\sqrt{n}})^2$

$\hspace{1.5cm}\Rightarrow s_{Q}(\mathcal{H})(q_{\max}+q_{\min})(\tau^2-\frac{1}{n})<\tau^2q^2_{\max}$, using Theorem $\ref{t2.13}$

$\hspace{1.5cm}\Rightarrow s_{Q}(\mathcal{H})<\frac{\tau^2}{(\tau^2-\frac{1}{n})}\frac{q^2_{\max}}{(q_{\max}+q_{\min})}$

Hence using Theorem $\ref{t2.13}$ and Theorem $\ref{thm2.12}$ we get, $s_{Q}(\mathcal{H})<\frac{2\tau^2d_{\max}}{(\tau^2-\frac{1}{n})}$.\qe

\subsection{Signless Laplacian spread of hypergraphs in terms of first Zagreb index}

In this section, we shall supply bounds for the signless Laplacian spread of a $k$-uniform hypergraphs in terms of first Zagreb index. Let $Z_1(\mathcal{H})=\sum\limits_{i\in V}d^2_{i}$ be call the first Zagreb index of a hypergraph $\mathcal{H}$, where $d_i$ are degrees of $\mathcal{H}$, $i=1,2,\ldots,n$. To find the bounds for the signless Laplacian spread of a $k$-uniform hypergraphs, we use to denote $\alpha=\sum\limits_{i\in V}\sum\limits_{i\sim j}d^2_{ij}$.

\begin{theorem}\label{t3.5}
Let $\mathcal{H}$ be a $k$-uniform hypergraph with $n$ vertices and $m(\geq 1)$ hyperedges. Let $d_{\min}$ be the maximum degree of $\mathcal{H}$. Then 
\begin{center}
   $s_{Q}(\mathcal{H})\geq2d_{\min}-\sqrt{\frac{Z_{1}(\mathcal{H})+\frac{\alpha}{(k-1)^2}-(2d_{\min})^2}{n-1}}$,
\end{center}
Equality holds if and only if the form signless Laplacian matrix is $Q(\mathcal{H})= t((n-1)I_n+J_n)$, where $J_n$ is the $n\times n$ matrix with each entry $1$, $I_n$ is the identity matrix of order $n$ and $t$ is a positive real number.
\end{theorem}

\pr Since 
\begin{align*}
Z_{1}(\mathcal{H})+\frac{\alpha}{(k-1)^2}&=\sum\limits_{i=1}^{n}q_{i}(\mathcal{H})^2\\
&\geq q_{\max}(\mathcal{H})^2+(n-1)q_{\min}(\mathcal{H})^2\\
q_{\min}(\mathcal{H})^2&\leq\frac{Z_{1}(\mathcal{H})+\frac{\alpha}{(k-1)^2}-q_{\max}(\mathcal{H})^2}{n-1}
\end{align*}

Therefore $q_{\min}(\mathcal{H})\leq\sqrt{\frac{Z_{1}(\mathcal{H})+\frac{\alpha}{(k-1)^2}-q_{\max}(\mathcal{H})^2}{n-1}}$

We know $s_{Q}(\mathcal{H})=q_{\max}(\mathcal{H})-q_{\min}(\mathcal{H}).$ Then we have
\begin{center}
$s_{Q}(\mathcal{H})\geq q_{\max}(\mathcal{H})-\sqrt{\frac{Z_{1}(\mathcal{H})+\frac{\alpha}{(k-1)^2}-q_{\max}(\mathcal{H})^2}{n-1}}$
\end{center}
Using Theorem $\ref{thm2.12}$ we have,\begin{center}
    $s_{Q}(\mathcal{H})\geq2d_{\min}-\sqrt{\frac{Z_{1}(\mathcal{H})+\frac{\alpha}{(k-1)^2}-(2d_{\min})^2}{n-1}}.$
\end{center}

The proof the equality part is same as the proof of Theorem \ref{thm:2.8}.
\qe

\subsection{Signless Laplacian spread of hypergraphs in terms of chromatic number}
In \cite{lima2011}, the authors proved that for a graph $G$ of order $n$, $s_Q(G)\geq\chi(G)$ and equality holds if and only if $G$ is a complete graph. This result is not true in general for $k(\geq 3)$-uniform hypergraphs. For example, consider the hypergraph $\mathcal{H}=(V,E)$ where $V=\{1,2,3\}$ and $E=\big\{\{1,2,3\}\big\}.$ For this hypergraph $s_Q(\mathcal{H})=1.5$ and $\chi(\mathcal{H})=3.$ Immediately one can think that the bound can be $s_Q(G)\geq\frac{\chi(G)}{k-1}$ or 
$s_Q(G)\geq (k-1)\chi(G).$ But this not true in general. For example consider the hypergraph $\mathcal{H}=(V,E)$ where $V=\{1,2,3\}$ and $E=\big\{\{1,2,3\},\{1,4,5\}\big\}.$ For this hypergraph above two bounds are not true. 

In this section we obtain an upper bound and a lower bound for the signless Laplacian spread in terms of strong chromatic number. To proceed further we need the following lemma.

\begin{lemma}{\rm\cite{izumino1998}}\label{lem3.1}
If $0\leq m_{1}\leq a_i\leq M_{1}$ and $0\leq m_{2}\leq a_i\leq M_{2}$, $1\leq i\leq n$ are two real numbers, then 
\begin{center}
$\sum\limits_{i=1}^{n}a^2_i\sum\limits_{i=1}^{n}b^2_i-(\sum\limits_{i}^{n}a_i b_i)^2\leq\frac{n^2}{4}(M_1M_2-m_1m_2)^2,$
\end{center}
\end{lemma}

\begin{lemma}{\rm\cite{saha2022}}\label{l3.13}
Let $\mathcal{H}$ be a $k$-uniform hypergraph. Then $\chi(\mathcal{H})\leq 1+(k-1)d_{\max}$, where $d_{\max}$ is the maximum degree of $\mathcal{H}$.
\end{lemma}

\begin{theorem}
Let $\mathcal{H}$ be a connected $k$-uniform hypergraph with $n(\geq 2)$ vertices. Then \begin{center}
    $\frac{\chi(\mathcal{H})\sqrt{n^2-1}}{n(1+(k-1)d_{\max})}<s_Q(\mathcal{H})<\frac{4n\chi(\mathcal{H})d_{\max}}{k\sqrt{n^2-1}}$,
\end{center}
where $d_{\max}$ and $\chi(\mathcal{H})$ are the maximum degree and strong chromatic number of $\mathcal{H}$, respectively.

\end{theorem}

\pr Let $a_i=\frac{s_Q(\mathcal{H})}{\chi(\mathcal{H})}i$, $b_i=1$, for $i=1,2,\ldots,n$. By using Theorem $\ref{thm2.12}$ and Lemma $\ref{l3.13}$, we have 
\begin{center}
$\frac{1}{1+(k-1)d_{\max}}\leq a_{i}\leq\frac{2nd_{\max}}{k}$, for $i=1,2,\ldots,n$.
\end{center} 

Let $M_{1}=\frac{2nd_{\max}}{k}$ and $m_{1}=\frac{1}{1+(k-1)d_{\max}}$.
Using Lemma $\ref{lem3.1}$ we get,
\begin{align*}
\sum\limits_{i=1}^{n}\bigg(\frac{s_Q(\mathcal{H})}{\chi(\mathcal{H})}i\bigg)^2\sum\limits_{i=1}^{n}1^2-\bigg(\sum\limits_{i}^{n}\frac{s_Q(\mathcal{H})}{\chi(\mathcal{H})}i\bigg)^2&\leq\frac{n^2}{4}\bigg(\frac{2nd_{\max}}{k}-\frac{1}{1+(k-1)d_{\max}}\bigg)^2\\
\Rightarrow \bigg(\frac{s_Q(\mathcal{H})}{\chi(\mathcal{H})}\bigg)^2\bigg(n\sum\limits_{i=1}^{n}i^2-(\sum\limits_{i}^{n}i)^2\bigg)&\leq\frac{n^2}{4}\bigg(\frac{2nd_{\max}}{k}-\frac{1}{1+(k-1)d_{\max}}\bigg)^2\\
\Rightarrow \bigg(\frac{s_Q(\mathcal{H})}{\chi(\mathcal{H})}\bigg)^2\frac{n^2(n^2-1)}{12}&\leq\frac{n^2}{4}\bigg(\frac{2nd_{\max}(1+(k-1)d_{\max})-k}{k(1+(k-1)d_{\max})}\bigg)^2\\
&<\frac{n^2}{4}\bigg(\frac{2nd_{\max}(1+(k-1)d_{\max})}{k(1+(k-1)d_{\max})}\bigg)^2\\
\Rightarrow \frac{s_{Q}(\mathcal{H})}{\chi(\mathcal{H})}&<\frac{4n^2d_{\max}}{4k}\frac{\sqrt{12}}{n\sqrt{n^2-1}}
\end{align*}
Hence $s_Q(\mathcal{H})<\frac{4n\chi(\mathcal{H})d_{\max}}{k\sqrt{n^2-1}}$.

Again let $a_i=\frac{\chi(\mathcal{H})}{s_Q(\mathcal{H})}i$, $b_i=1$, for $i=1,2,\ldots,n$. By using Theorem $\ref{thm2.12}$ and Lemma $\ref{l3.13}$, $\frac{k}{2d_{\max}}\leq a_{i}\leq n(1+(k-1)d_{\max})$.
Let $M_{1}=n(1+(k-1)d_{\max})$ and $m_{1}=\frac{k}{2d_{\max}}$.
Again using Lemma \ref{lem3.1} we get, 
\begin{align*}
\sum\limits_{i=1}^{n}\bigg(\frac{\chi(\mathcal{H})}{s_Q(\mathcal{H})}i\bigg)^2\sum\limits_{i=1}^{n}1^2-\bigg(\sum\limits_{i}^{n}\frac{\chi(\mathcal{H})}{s_Q(\mathcal{H})}i\bigg)^2&\leq\frac{n^2}{4}\bigg(n(1+(k-1)d_{\max})-\frac{k}{2d_{\max}}\bigg)^2\\
\Rightarrow \bigg(\frac{\chi(\mathcal{H})}{s_Q(\mathcal{H})}\bigg)^2\frac{n^2(n^2-1)}{12}&\leq\frac{n^2}{4}\bigg(\frac{nd_{\max}(1+(k-1)d_{\max})-k}{2d_{\max}}\bigg)^2\\
&<\frac{n^4(1+(k-1)d_{\max})^2}{16}\\
\Rightarrow \bigg(\frac{\chi(\mathcal{H})}{s_Q(\mathcal{H})}\bigg)^2&<\frac{n^4(1+(k-1)d_{\max})^2}{16}\frac{12}{n^2(n^2-1)}
\end{align*}
Hence $s_Q(\mathcal{H})>\frac{\chi(\mathcal{H})\sqrt{n^2-1}}{n(1+(k-1)d_{\max})}$.\qe

\end{document}